\documentclass[10pt]{article}
\usepackage[margin=1in]{geometry}

\usepackage[T1]{fontenc}
\usepackage{lmodern}
\usepackage{todonotes}
\usepackage{ulem}
\usepackage{bm}
\usepackage[colorlinks=true,linkcolor=red,citecolor=blue]{hyperref}%
\usepackage{cite}
\usepackage{url}            
\usepackage{booktabs}       
\usepackage{amsfonts}    
\usepackage{nicefrac}       
\usepackage{microtype} 
\usepackage{wrapfig}
\usepackage{enumerate}
\usepackage{graphicx}
\usepackage{amsmath,amssymb}
\usepackage{bigstrut}
\usepackage{float}
\usepackage{caption}
\usepackage{algorithm}
\usepackage{mathtools}
\usepackage{subfigure}
\usepackage{booktabs}
\usepackage{authblk}
\usepackage{algpseudocode}
\usepackage{dsfont}
\usepackage{bbm}

\newtheorem{theorem}{Theorem}

\newtheorem{remark}{Remark}
\newtheorem{lemma}{Lemma}

\newtheorem{prop}{Proposition}

\newtheorem{asm}{Assumption}

\def\R{\mathbb{R}}
\def\ones{\mathbf{1}}
\def\xh{\hat{\bfx}}
\def\nt{\bfe}
\def\sm{\gamma}

\newcommand{\EE}[1]{\mathbb{E}\left[ #1 \right]}

\newcommand{\bfx}{\mathbf{x}}
\newcommand{\bfe}{\mathbf{e}}
\newcommand{\bfy}{\mathbf{y}}
\newcommand{\bfv}{\mathbf{v}}

\newcommand{\bfr}{\mathbf{r}}

\newcommand{\bfg}{\mathbf{g}}

\newcommand{\bfu}{\boldsymbol{u}}

\newcommand{\clg}{\mathcal{G}}

\newcommand{\cle}{\mathcal{E}}

\def\Pr{\mathop{\rm Pr}\nolimits}
\def\V{V_\bfr}
\def\rmin{\bfr_{\min}}
\def\Wt{\{W(t)\}}
\def\minW{\eta}

\newcommand{\comment}[1]{}

\newcommand{\F}{\mathcal{F}}
\newcommand{\azr}{\alpha_0}
\newcommand{\bzr}{\beta_0}

\newcommand{\md}{\middle| }
\DeclareMathOperator{\diag}{diag}

\newcommand{\mx}{\bar{\bfx}}

\newcommand{\lp}{\left(}
\newcommand{\rp}{\right)}
\newcommand{\lb}{\left[}
\newcommand{\rb}{\right]}

\newcommand{\lnr}{\left\|}
\newcommand{\rnr}{\right\|}
\newcommand{\lc}{\left\{}
\newcommand{\rc}{\right\}}
\newcommand{\dl}{\delta}

\newcommand{\BO}{\mathcal O}

\newcommand{\ER}[1]{E(#1)}

\newcommand{\cOne}{c_1}
\newcommand{\cTwo}{c_2}
\newcommand{\cThree}{c_3}
\newcommand{\cFour}{c_4}

\newcommand{\Nr}[1]{\lnr #1 \rnr_{\bfr}}

\def\R{\mathbb{R}}

\def\E{\mathbb{E}}

\def\G{G}

\newcommand{\norm}[1]{\left\|#1\right\|}

\newcommand{\cR}{\mathcal{R}}

\newcommand{\Dl}{D}
\newcommand{\tz}{t_0}
\newcommand{\kz}{k_0}
\newcommand{\tone}{t_1}
\newcommand{\ttwo}{t_2}

\newcommand{\cfone}{S_1}
\newcommand{\cftwo}{S_2}
\renewcommand{\k}{t}

\title{Almost Sure Convergence of Distributed Optimization with Imperfect Information Sharing}

\author{Hadi Reisizadeh, Anand Gokhale, Behrouz Touri, and Soheil Mohajer  
\thanks{ H.\ Reisizadeh (email: hadir@umn.edu) and S.\ Mohajer (email: soheil@umn.edu) are with the University of Minnesota, A.\ Gokhale (email:anand\_gokhale@ucsb.edu) is with the University of California Santa Barbara, and B.\ Touri (email: btouri@ucsd.edu) is with the University of California San Diego.}}

\date{}
\begin{document}

\maketitle

\begin{abstract} 
To design algorithms that reduce communication cost or meet rate constraints and are robust to communication noise, we study convex distributed optimization problems where a set of agents are interested in solving a separable optimization problem collaboratively with imperfect information sharing over time-varying networks.  We study the almost sure convergence of a two-time-scale decentralized gradient descent algorithm to reach the consensus on an optimizer of the objective loss function. One time scale fades out the imperfect incoming information from neighboring agents, and the second one adjusts the local loss functions' gradients. We show that under certain conditions on the connectivity of the underlying time-varying network and the time-scale sequences, the dynamics converge almost surely to an optimal point supported in the optimizer set of the loss function.
\end{abstract}
 
\section{Introduction}

In recent years, with the rapid growth of areas such as big data and machine learning, there has been a surge of interest in  studying multi-agent networks. In many machine learning applications, it is impractical to implement the learning task in a
centralized fashion due to the decentralized nature of datasets. Multi-agent networked systems provide scalability to larger datasets and systems, data locality, ownership, and privacy, especially for modern computing services. These systems arise in various applications such as sensor networks~\cite{rabbat2004distributed}, multi-agent control~\cite{olshevsky2010efficient}, large-scale machine learning~\cite{tsianos2012consensus}, and power networks~\cite{ram2009distributed}. The task of learning a common objective function over a multi-agent network can be reduced to a distributed optimization problem.

In distributed optimization problems, a set of agents are interested in finding a minimizer of a separable function $f(\bfx) := \sum_{i=1}^{n} f_i(\bfx)$ such that each agent $i$ has access to the local and private loss function $f_i(\cdot)$ of this decomposable cost function. Therefore, the goal is to identify system dynamics that guarantee the asymptotic convergence of all agent states to a common state which minimizes the objective cost function $f(\cdot)$. Various methods have been proposed to solve distributed optimization problems in strongly convex~\cite{nedic2017achieving,saadatniaki2018optimization}, convex~\cite{nedic2009distributed,nedic2010constrained}, and non-convex settings~\cite{di2016next,tatarenko2016local,tatarenko2017non,hong2017prox}.  Under perfect information sharing assumption, a subgradient-push algorithm is proposed for (strongly) convex loss functions in~\cite{nedic2014distributed}. The almost sure and $L_2$-convergences of this method are shown for certain conditions on the connectivity, loss functions, and the step-size sequence.

Most of the works mentioned above on distributed optimization problems assume ideal  communication channels and perfect information sharing among agents. However, most communication channels are noisy, and exchanging real-valued vectors introduces a massive communication overhead. To mitigate this challenge, various compression approaches have been introduced~\cite{alistarh2017qsgd,wangni2018gradient} where the vectors are represented with a finite number of bits or only a certain number of the most significant coordinates of the vectors are selected. Then, the quantized/sparsified vectors are communicated over the network. Consequently, each agent receives an imperfect estimate of the intended messages from the neighboring agents. Various gradient descent algorithms with imperfect information sharing have been proposed~\cite{koloskova2019decentralized,koloskova2019decentralizedICLR,reisizadeh2021distributed}, showing the convergence rates in $L_2$ sense for the diverse set of problem setups. 
 
A related work~\cite{srivastava2011distributed} considers distributed \textit{constrained} optimization problems with noisy communication links. The work differs from our current study as they assumed \textit{i.i.d.}\ weight matrices with a symmetric expected weight matrix. In a recent work~\cite{doan2020convergence}, a two-time-scale gradient descent algorithm has been presented for (strongly) convex cost functions over a convex compact set. Assuming a \textit{fixed} topology for the underlying network, the uniform weighting of the local cost functions, and a \textit{specific} scheme for the  lossy sharing of information, it is shown that under certain conditions, the agents' states converge to the optimal solution of the problem almost surely. Considering the averaging-based distributed optimization over  random networks with possible dependence on the past and under certain conditions, the almost sure convergence to an optimal point is presented in~\cite{aghajan2020distributed}. 

In this paper, we study the distributed \textit{convex} optimization problems over \textit{time-varying} networks with \textit{imperfect} information sharing. 
We consider the two-time-scale gradient descent method studied in~\cite{reisizadeh2021distributed,reisizadeh2022dimix} to solve the optimization problem. One time-scale adjusts the (imperfect) incoming information from the neighboring agents, and one time-scale controls the local cost functions' gradients. It is shown that with a proper choice of parameters, the proposed algorithm reaches the global optimal point for \textit{strongly} convex loss functions at a rate of $\BO(T^{-1/2})$ and achieves a convergence rate of $\BO(T^{-1/3+\epsilon})$ with any $\epsilon>0$ for non-convex cost function in $L_2$ sense. 
Here, we identify the sufficient conditions on the step-sizes sequences for the almost sure convergence of the agent’s states to an optimal solution for the class of convex cost functions.

The paper is organized as follows. We conclude this section by discussing the  notations used in the paper. We formulate the main problem and state the relevant underlying  assumptions in Section~\ref{sec:ProblemFormulation}. In Section~\ref{sec:main}, we present our main results. To corroborate our theoretical analysis, we present simulation results in Section~\ref{sec:sim}. We discuss some preliminary results which are required to prove the main results in Section~\ref{sec:prelims}. The proofs of the preliminary lemmas are presented in Appendix. Then, we present the proof of the main results in Section~\ref{proof:main_result} and Section~\ref{sec:proof-prop}.
Finally, we conclude the paper and discuss some possible directions for future works in Section~\ref{sec:conclusion}. 

{\bf Notation}. We denote the set of integers $\{1,\ldots,n\}$ by $[n]$. In this paper, we consider $n$ agents that are minimizing a function in $\R^d$. We assume that all local objective functions are acting on \textbf{row} vectors in $\R^{1\times d} = \R^d$, and thus we view vectors in $\R^d$ as row vectors. Note that the rest of the vectors, i.e., the vectors in $\R^{n\times 1}=\R^n$, are assumed to be column vectors. We denote the $L_2$ norm of a vector $\bfx\in \R^d$ by $\norm{\bfx}$. For a matrix $A \in \R^{n\times d}$, and a strictly positive stochastic vector $\bfr \in \R^n$, we define the $\bfr-$norm of $A$ by $\norm{A}^2_\bfr = \sum_{i=1}^{n} r_i \norm{A_i}^2_2$, where $A_i$ denotes the $i$-th row of $A$. We denote the Frobenius norm for a matrix $A\in \R^{n\times d}$ by $\norm{A}^2 = \sum_{i=1}^{n}\sum_{j=1}^{d} |A_{ij}|^2$. The difference operator $\Delta$ on a real-valued sequence $\{a(t)\}$ is defined as $\Delta a(t):= a(t+1)-a(t)$ for every $t\geq 1$. 


\section{Problem Formulation}\label{sec:ProblemFormulation}
In this section, we first formulate the problem of interest. Then, we present the underlying assumptions on the information exchange, the network connectivity, the properties of the cost functions, and the time-scales sequences. 

\subsection{Problem Statement}
Consider a set of $n\geq 2$ agents, that are connected through a time-varying network. Each agent $i \in [n]$ has access to a local cost function $f_i: \R^d \rightarrow \R$. The goal of this work is to solve the optimization problem which is given by
\begin{align}\label{eqn:MainProblem}
    \!\!\min_{\bfx_1, \bfx_2, \ldots, \bfx_n \in \R^d} \sum_{i=1}^{n}r_i f_i(\bfx_i) \quad\textrm{subject to}\quad  \bfx_1 =\cdots= \bfx_n,\nonumber
\end{align}
where $\bfr = (r_1,r_2,\ldots,r_n)^T$ is a stochastic vector, i.e., $r_i\geq 0$ and $\sum_{i=1}^{n} r_i=1$.

At each iteration $t\geq 1$,  we represent the topology of the network connecting the agents by a directed graph ${\clg(t) = ([n],\enspace \cle(t))}$, where the vertex set $[n]$ represents the agents, and the edge set $\cle(t) \subseteq [n]\times [n]$ represents the connectivity pattern of the agents at time $t$, i.e.,  the edge ${(i,j) \in \cle(t)}$ denotes a directed edge from agent $i$ to agent $j$. We assume that at each time $t\geq 1$, each agent $i\in[n]$ can only send a message to its out-neighbours, i.e., the set of all agents $j$ such that $(i,j)\in \cle(t)$. We also assume that $\clg(t)$ satisfies certain connectivity conditions, which are discussed in Assumption~\ref{Assumption:Graph}.

In this work, the communication between the agents is assumed to be imperfect. We adapt the general framework of the noisy sharing of information introduced in \cite{reisizadeh2021distributed} as described below. Given the states $\bfx_i(t)$ of agents $i\in [n]$ at time $t$, we assume that each agent has access to an imperfect weighted average of its in-neighbours states, denoted by  $\xh_i(t)$ given by 
\[{\xh_i(t) = \sum_{j=1}^{n}W_{ij}(t)\bfx_j(t) + \bfe_i(t)},\]
where $W(t) = [W_{ij}(t)]$ is a row-stochastic matrix and $\bfe_i(t)$ is a random noise vector in $\R^d$. Note that the matrix $W(t)$ is consistent with the underlying graph $\clg(t)$. More precisely, we have  $W_{ij}(t)> 0$ if and only if $(j,i) \in \cle(t)$, where $\cle(t)$ is the edge set of the graph $\clg(t)$. 

Regarding the local cost function, we assume that agent $i \in [n]$ has access to a subgradient $\bfg_i(\bfx_i(t))$ of the local cost function $f_i(\cdot)$ at each local decision variable $\bfx_i(t)$ at time $t$.
Inspired by~\cite{reisizadeh2021distributed}, we present the update rule in this work as
\begin{equation} \label{eqn:UpdateLawXHat}
    \bfx_i(t+1) = (1 - \beta(t))\bfx_i(t) + \beta(t) \xh_i(t) - \alpha(t)\bfg_i(\bfx_i(t)),
\end{equation}
where $\{\beta(t)\}$ and $\{\alpha(t)\}$ are the sequences of step-sizes of the algorithm. This work identifies \textit{sufficient} conditions on the sequences of step-sizes for almost sure convergence of the dynamics~\eqref{eqn:UpdateLawXHat} to an optimal point ${\bfx^\star \in \mathcal{X}^\star = \arg \min_{\bfx\in \R^d}\sum_{i=1}^n f_i(\bfx)}$. 
For simplicity of notation, let
 \begin{align}
        X(t) :=\begin{bmatrix} 
        \!\bfx_1(t)\! \\ \vdots \\ \!\bfx_n(t)\!
        \end{bmatrix}, 
        E(t):= \begin{bmatrix}
        \!\nt_1(t)\! \\ \vdots\\ \!\nt_n(t)\!\end{bmatrix},
        \G(t) := \begin{bmatrix} \!\bfg_1(\bfx_1(t))\!\\ \vdots \\ \!\bfg_n(\bfx_n(t))\!\end{bmatrix}.\nonumber
\end{align}
Using these matrices, we can write the update rule \eqref{eqn:UpdateLawXHat} in the form of a linear time-varying system given by 
\begin{align}\label{eqn:UpdateLawMatrixForm}
    X(t+1) &=A(t)X(t) + U(t),
\end{align}
where 
\[A(t) :=(1 - \beta(t))I + \beta(t)W(t)\]
and 
\[U(t):=\beta(t) E(t) - \alpha(t) \G(t).\]
We also define ${\Phi(t\hspace{-1pt}:\hspace{-1pt}s):=A(t\hspace{-1pt}-\hspace{-1pt}1)\cdots A(s\hspace{-1pt}+\hspace{-1pt}1)}$  for $t> s$, with  $\Phi(t\hspace{-1pt}:\hspace{-1pt}t-1)=I$.
\subsection{Assumptions}
To proceed with our main result, we need to make certain assumptions regarding the noise vectors $\{\bfe_i(t)\}$, the weight matrix $\{W(t)\}$, the local cost functions, and the sequences of step-sizes. 
\begin{asm}[Noise Sequence Assumptions]\label{Assumption:Noise}
We assume that the noise sequence $\{\bfe_i(t)\}$ satisfies 
\begin{align*}
    \EE{\bfe_i(t)\md\F_t} = 0, \enspace \textrm{and} \enspace
    \EE{\norm{\bfe_i(t)}^2\md\F_t} \leq \gamma,
\end{align*}
for some $\gamma > 0$, all $i \in [n]$, and all $t\geq 1$. Here, $\{\F_t\}$ is the natural filtration of the random process $\{X(t)\}$.
\end{asm}

\begin{asm}[Connectivity Assumptions]\label{Assumption:Graph}
We assume that the sequence $\{\clg(t)\}$ of underlying graphs, and its associated weight matrix sequence $\{W(t)\}$ satisfy the following properties.
\begin{enumerate}[(a)]
    \item $W(t)$ is non-negative, $W(t)\ones = \ones$, and ${\bfr^T W(t) = \bfr^T}$ for all $t\geq 1$, where $\ones$ is the all-one vector, and $\bfr>0$ is the given stochastic weight vector.
    \item Each nonzero element of $W(t)$ is bounded away from zero, i.e., there exists some $\eta >0$, such that if ${W_{ij}(t) > 0}$ for $i,j \in [n]$ and $t \geq 1$, then $W_{ij}(t) \geq \eta$.
    \item There exists an  integer $B \geq 1$ such that the graph $\left([n],\bigcup_{k=t+1}^{t+B}\cle(k) \right)$ is strongly connected for all $t \geq 1$, where $\cle(k)$ is the edge set of $\clg(k)$.
\end{enumerate}
\end{asm}

\begin{asm}[Objective Function Assumptions]\label{Assumption:Function}
We assume that objective functions $f_i$ satisfy the following properties.
\begin{enumerate}[(a)]
    \item $f_i$ is convex for all $i \in [n]$.
    \item  The optimizer set $\mathcal{X}^\star :=\arg \min_{\bfx\in \mathbb{R}^d} \sum_{i=1}^n r_i f_i(\bfx) $ is non-empty.
    \item  Each $f_i$ has bounded subgradients, i.e., there exists ${L > 0}$ such that $\norm{\bfg_i} < L$ for all subgradients $\bfg_i$ of $f_i(\bfx)$ at every $\bfx\in \R^d$.
    This also implies that each $f_i(\cdot)$ is $L$-Lipschitz continuous, i.e., 
    \[|f_i(\bfx) - f_i(\bfy)|< L\norm{\bfx - \bfy},\]
    for all $\bfx,\bfy\in \R^d$.
\end{enumerate} 
\end{asm} 
\begin{asm}[Step-size Sequences Assumptions]\label{asm:step-size}
For the non-increasing step-size sequences $\{\alpha(t)\}$ and $\{\beta(t)\}$ where $\{\beta(t)\}$ take values in $[0,1]$, we assume that 
\begin{enumerate}[(a)]
\item \label{asm:4-a} $\sum_{t=1}^\infty \alpha(t) = \infty$,
\item \label{asm:4-a2-b2} $\sum_{t=1}^\infty \alpha^2(t) < \infty$, $\sum_{t=1}^\infty \beta^2(t) < \infty$, and
\item \label{asm:4-a2/b}$\sum_{t=1}^\infty \frac{\alpha^2(t)}{\beta(t)}<\infty$.
\end{enumerate}
Also, there exists some $t_0\geq 1$ such that for every $t\geq t_0$  
\begin{enumerate}[(a)]
\setcounter{enumi}{3}
\item \label{asm:4-Db} $-\Delta \beta(t) \leq \cOne \beta^2(t)$, 
\item \label{asm:4-Da} $-\Delta \alpha(t)\leq \cTwo \alpha(t)\beta(t)$,
\end{enumerate} 
for some positive constants $\cOne\!<\!\frac{\lambda}{2}$ and $\cTwo\!<\!\frac{\lambda}{4}$, where $\lambda:=\frac{\eta\rmin}{2Bn^2}<1$ and $\rmin:= \min_{i\in [n]} \{r_i\}\leq 1$.
\end{asm}

\begin{remark}\label{rem:asm:step-size}
Note that Assumptions~\ref{asm:step-size}-\eqref{asm:4-a},~\eqref{asm:4-a2-b2}, and~\eqref{asm:4-a2/b}  imply $\sum_{t=1}^\infty\beta(t)=\infty$ as if $\sum_{t=1}^{\infty} \beta(t)<\infty$, using the Cauchy-Schwarz inequality we get
\begin{align*}
     \sum_{t=1}^{\infty} \alpha(t) \leq \lp \sum_{t=1}^{\infty} \frac{\alpha^2(t)}{\beta(t)}\rp^{\frac{1}{2}}\lp\sum_{t=1}^{\infty} \beta(t)\rp^{\frac{1}{2}} <\infty,
\end{align*}
which is a contradiction with Assumption~\ref{asm:step-size}-\eqref{asm:4-a}.
Similarly, we can write
\begin{align}\label{eq:sum_a_b_hf}
    \sum_{t=1}^{\infty} \alpha(t)\beta^{\frac{1}{2}}(t) \leq \lp \sum_{t=1}^{\infty} \frac{\alpha^2(t)}{\beta(t)}\rp^{\frac{1}{2}}\lp\sum_{t=1}^{\infty} \beta^2(t)\rp^{\frac{1}{2}} <\infty,
\end{align}
where the second inequality follows from Assumptions~\ref{asm:step-size}-\eqref{asm:4-a2-b2} and~\eqref{asm:4-a2/b}.
\end{remark}
\begin{remark}
    Note that unlike Assumption~\ref{asm:step-size}-\eqref{asm:4-a}-\eqref{asm:4-a2/b} that do not depend on the dynamics parameters, Assumption~\ref{asm:step-size}-\eqref{asm:4-Db}-\eqref{asm:4-Da} depend on those parameters. However, we will show in Section~\ref{sec:proof-prop} that Assumption~\ref{asm:step-size}-\eqref{asm:4-Db}-\eqref{asm:4-Da} will be satisfied for sufficiently large $t$, regardless of the dynamic parameters.
\end{remark}
\section{Main Results}\label{sec:main}
In this section, we provide the main results of the paper. First, we present sufficient conditions for the sequences $\{\beta(t)\}$ and $\{\alpha(t)\}$ for the almost sure convergence to an optimal point for the agents acting under the dynamics~\eqref{eqn:UpdateLawXHat}. Then, for step-sizes of the form $\alpha(t) = \frac{\azr}{t^\nu}$ and $\beta(t) = \frac{\bzr}{t^\mu}$, we provide the region of $(\mu,\nu)$ for  which  the almost sure convergence is guaranteed. 

\begin{theorem}\label{thm:almost_sure}
If Assumptions~\ref{Assumption:Noise}-\ref{asm:step-size} are satisfied, then, for the dynamics~\eqref{eqn:UpdateLawXHat}, for all $i\in [n]$ we have ${\lim_{t\rightarrow \infty} \bfx_i(t) = \tilde{\bfx}}$ almost surely, where $\tilde{\bfx}$ is an optimal point in the set of optimal solutions $\mathcal{X}^\star$.
\end{theorem}
The proof of Theorem~\ref{thm:almost_sure} is provided in Section~\ref{proof:main_result}. The implication of the above result for the practical step-sizes $\alpha(t) = \frac{\azr}{t^\nu}$ and $\beta(t) = \frac{\bzr}{t^\mu}$ as as follows. 
\begin{prop}\label{prop:mu_nu}
Let  Assumptions~\ref{Assumption:Noise}-\ref{Assumption:Function} hold. Then, for every $i\in [n]$, the dynamics~\eqref{eqn:UpdateLawXHat} with step-sizes ${\alpha(t) = \frac{\alpha_0}{t^\nu}}$ and ${\beta(t) = \frac{\beta_0}{t^\mu}}$ converges almost surely, i.e., we have  ${\lim_{t\rightarrow \infty} \bfx_i(t) = \tilde{\bfx}}$ almost surely for some optimal point $\tilde{\bfx}\in\mathcal{X}^\star$,
provided that $\bzr\leq 1$, ${\frac{1}{2}<\mu \leq 1}$, and ${\frac{1}{2}(1\!+\!\mu)< \nu \leq 1}$.
\end{prop}
The proof of Proposition~\ref{prop:mu_nu} is provided in Section~\ref{sec:proof-prop}.
\begin{remark}
Proposition~\ref{prop:mu_nu} identifies sufficient conditions for $(\mu,\nu)$ such that the dynamics~\eqref{eqn:UpdateLawXHat} converges almost surely to the optimal set of a \textit{convex} objective function over \textit{time-varying} networks, when utilizing step-sizes of the form $\alpha(t) = \frac{\azr}{t^\nu}$ and $\beta(t) = \frac{\bzr}{t^\mu}$. The same dynamic is studied  for \textit{constrained} optimization problems with \textit{i.i.d.}\ weight matrices with symmetric expected weight,  and \textit{independent} noisy communication links in~\cite{srivastava2011distributed}, where 
interestingly, the same $(\mu,\nu)$-region for the convergence of the dynamic is obtained.   
\end{remark}

\begin{figure}[h]
	\centering
	\includegraphics[width=0.35\textwidth]{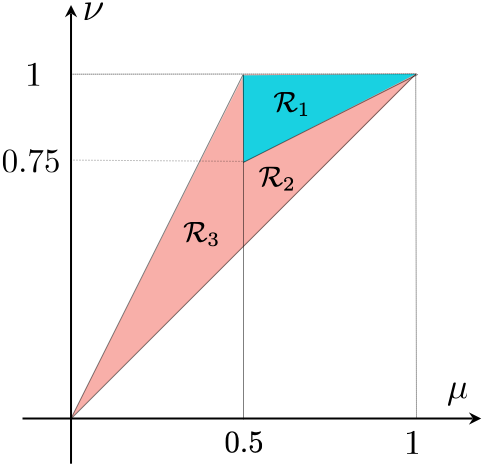}
	\caption{$\cR_1$ is the $(\mu,\nu)$-region for the almost sure convergence of the dynamic when applied on  strongly convex functions. The dynamic converges in the $L_2$-sense in $  \cR_1\cup \cR_2 \cup \cR_3$, when the objective functions are convex.} 
	\label{fig:Lu}
\end{figure}

\begin{remark}
In an earlier work~\cite{reisizadeh2021distributed}, the $L_2$-convergence of a similar dynamic is studied for \textit{strongly convex} objective functions.  Figure~\ref{fig:Lu} compares the the region of $(\mu,\nu)$ to guarantee the $L_2$-convergence~\cite{reisizadeh2021distributed} which is $\cR_1\cup\cR_2\cup\cR_3$ and the region for the almost sure convergence (this work), i.e., $\cR_1$. Interestingly, the optimal parameters $\mu^\star \!=\! 0.75$ and $\nu^\star \!=\! 1$ that lead to the fastest convergence in $L_2$ sense for strongly convex loss functions also guarantee the almost sure convergence.
\end{remark} 
Note that both regions only characterize  \textit{sufficient} conditions for two types of convergence criteria under different function properties, and hence, not comparable. For example, if for $L_2$ convergence, we relax the class of strongly convex functions to general convex functions, we can show that it is \textit{necessary} to have $\mu > \frac{1}{2}$. In other words, we cannot have $L_2$ convergence in region $\cR_3$ for the general class of (not necessarily strong) convex functions. To see this, consider the convex functions $f_i(\bfx)=0$ for $i\in [n]$ with the set of optimizers $\mathcal{X}^\star=\R^d$, and  zero-mean \textit{i.i.d.}\ noise sequences with variance $\gamma$, which satisfy Assumption~\ref{Assumption:Noise}. Then, multiplying both sides of~\eqref{eqn:UpdateLawMatrixForm} by $\bfr^T$, we get $\mx(k+1) = \mx(k) + \beta(k) \bfr^T \ER{k}$. 
Therefore, we can write
\begin{align}\label{eq:cons-beta}
    \mathbb{E}\left[\|\mx(k\!+\!1)\|^2\right] &\!=\! \EE{\EE{\|\mx(k\!+\!1)\|^2}\md \F_{k}}\nonumber\\
     &\!=\! \EE{\|\mx(k)\|^2} \!+\! \beta(k)\langle \bfr^T \EE{ \ER{k}\md\F_k}\!,\mx(k) \rangle + \beta^2(k)\EE{\left\|\bfr^T \ER{k}\right\|^2 \md\F_k}\nonumber \\
     &=\! \EE{\|\mx(k)\|^2} + \gamma \beta^2(k),
\end{align}
where the last equality is due to Assumption~\ref{Assumption:Noise} (see~\eqref{eq:quant-noise} for more details). Summing up~\eqref{eq:cons-beta} over $k$,  we arrive at
\begin{align*}
    \lim_{t\rightarrow \infty} \EE{\|\mx(t)\|^2} = \|\mx(1)\|^2 + \gamma \sum_{k=1}^{\infty} \beta^2(k).
\end{align*}
If $\bfx_i(t)$ converges to some $\tilde{\bfx}$ in $L_2$ for all $i\in[n]$, then we have $\lim_{t\rightarrow \infty} \EE{\|\mx(t)\|^2} = \EE{\|\tilde{\bfx}\|^2}<\infty$. This implies that  $\sum_{k=1}^{\infty} \beta^2(k)<\infty$, which means that we need to have $\mu>\frac{1}{2}$. In other words, the condition $\mu>\frac{1}{2}$ is necessary for the class of convex functions if $L_2$-convergence is desired. Further investigation is required to determine whether  region $\cR_2$ leads to a convergence.  
\section{Experimental Results}\label{sec:sim}

\begin{figure}[t]
    \centering
    \includegraphics[width=0.46\textwidth]{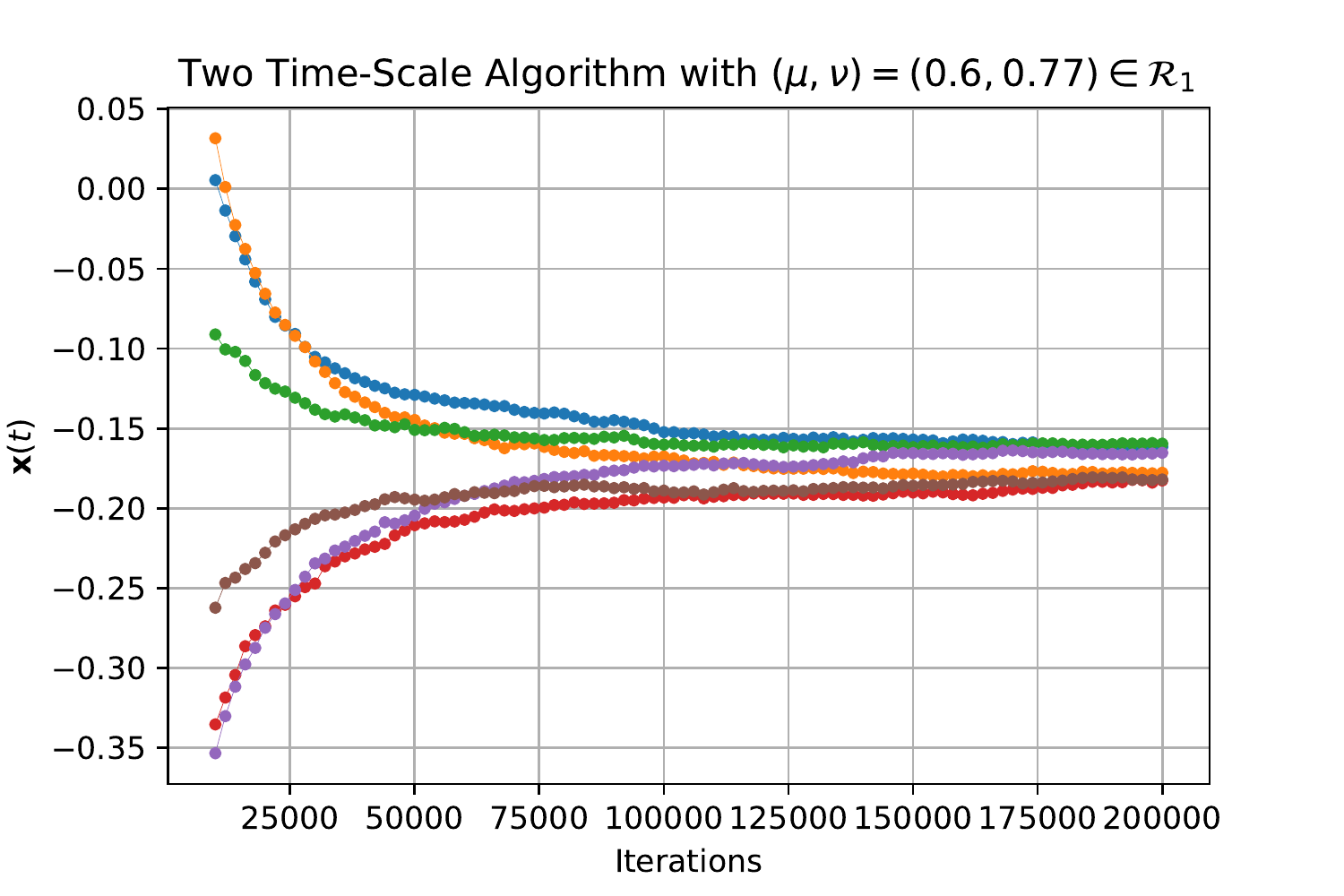}
    \vspace{-1mm}
\includegraphics[width=0.46\textwidth]{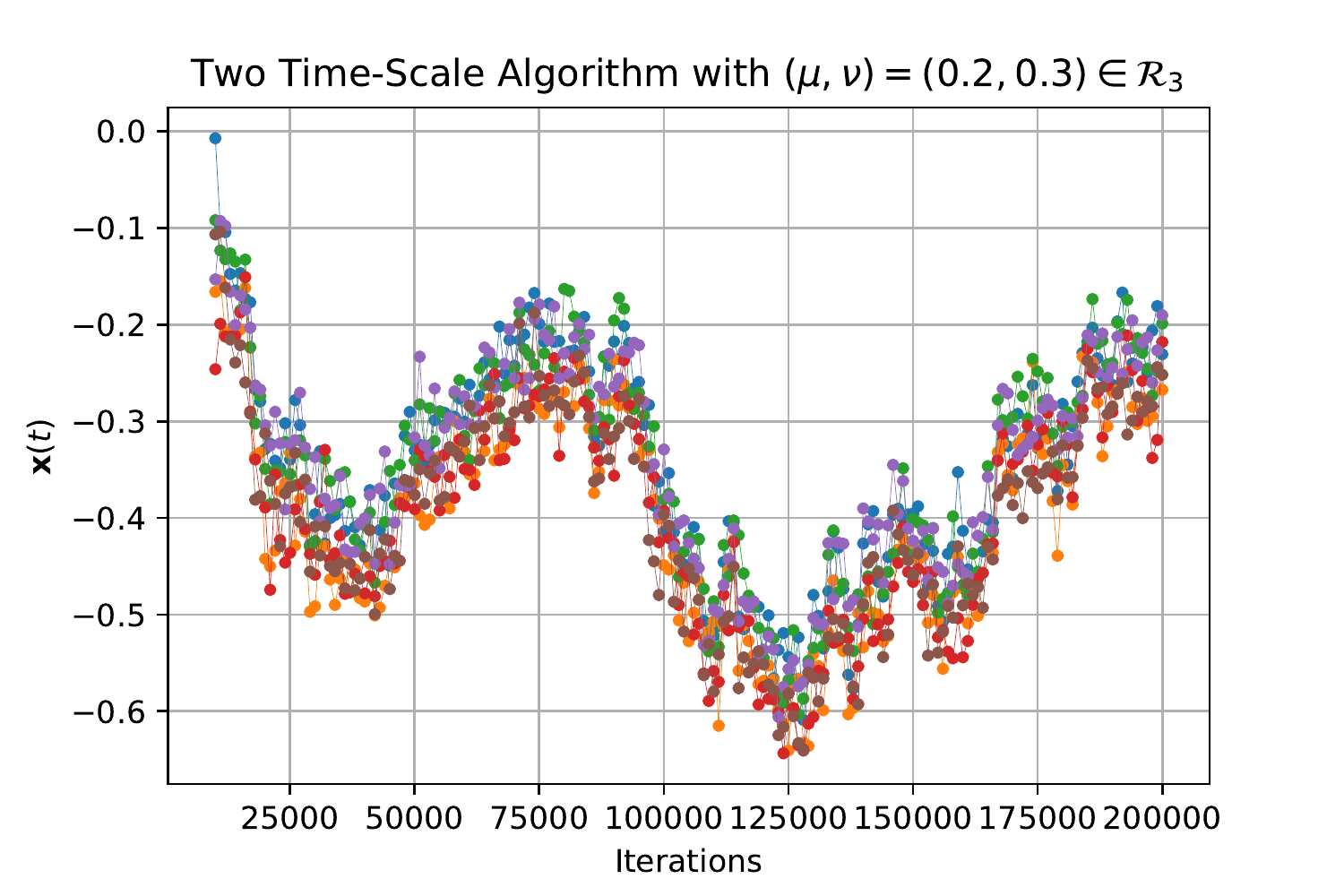}
     \vspace{-1mm}
     \includegraphics[width=0.46\textwidth]{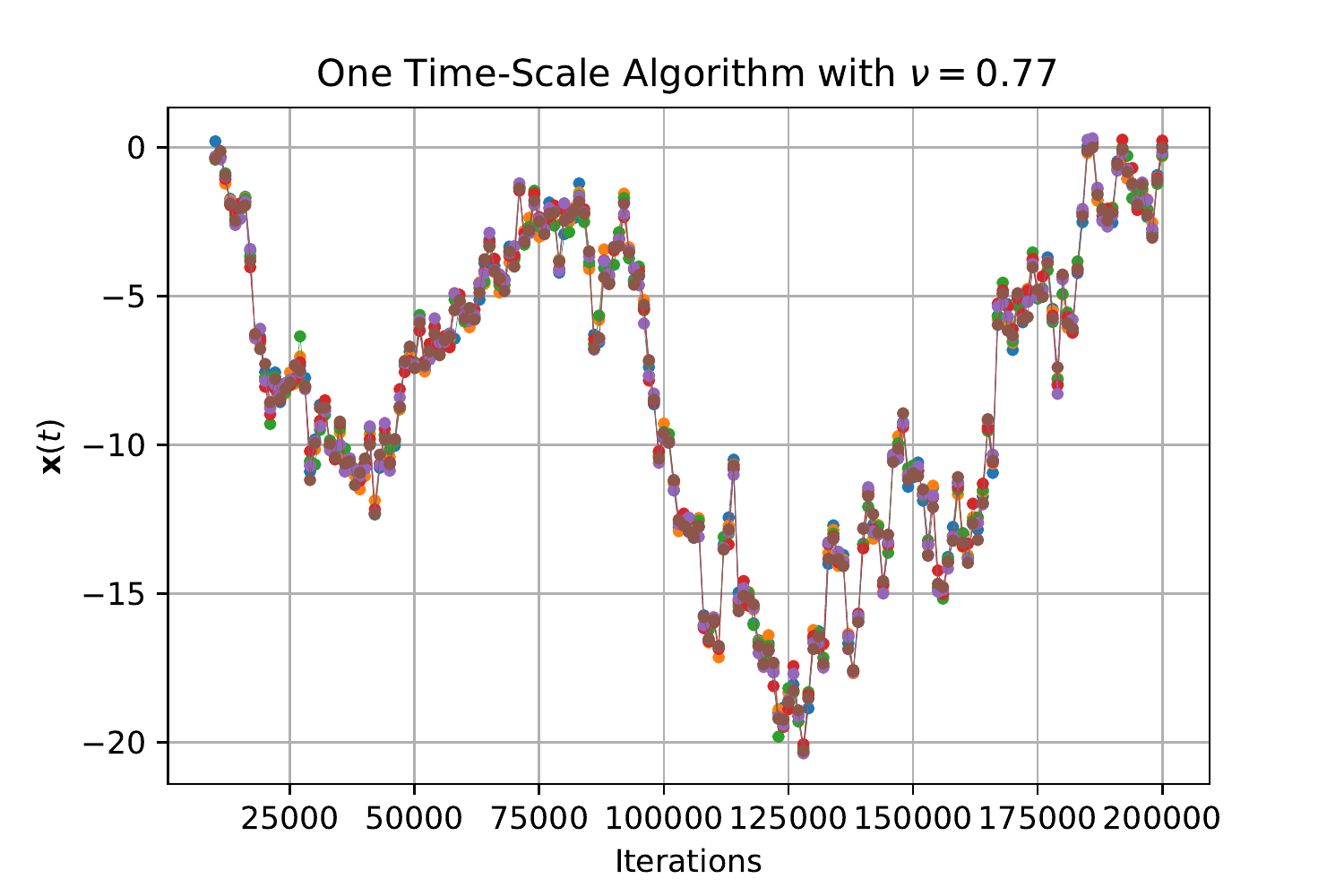}
     \vspace{-1mm}
\includegraphics[width=0.45\textwidth]{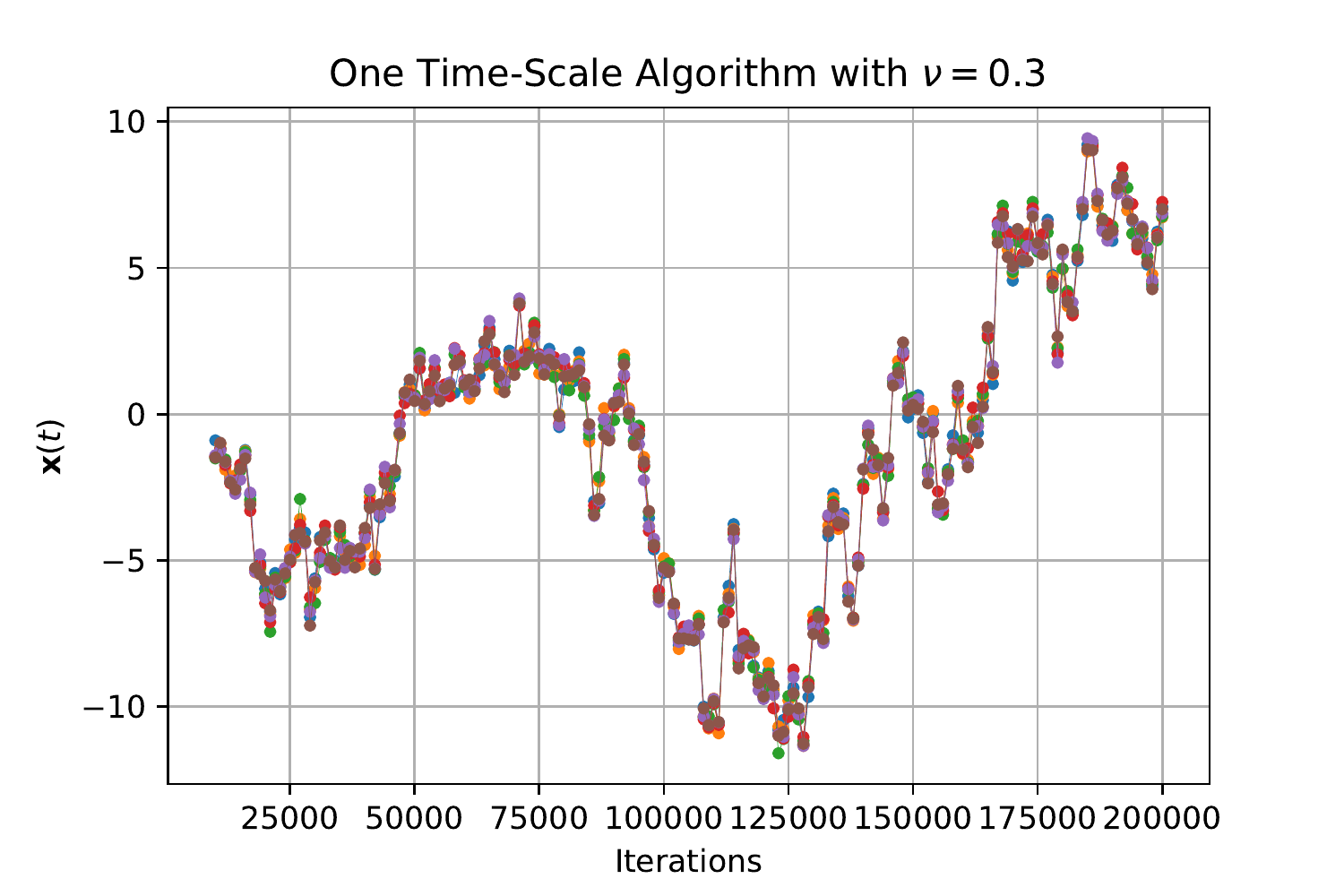}
     
\caption{Trajectory vs. Iterations: Two Time-Scale Algorithm with $(\mu,\nu) = (0.6,0.77)$ (top-left), $(0.2,0.3)$ (top-right), and One Time-Scale Algorithm with $\nu=0.77$ (bottom-left), $0.3$ (bottom-right).}
    \label{fig:sim}
    \vspace{-5mm}
\end{figure} 

In this section, we provide some numerical results to corroborate the derived theoretical analysis.

\textbf{Data and Experimental Setup.} We consider a time-varying network with $n=6$ agents, with lost functions given by $f_i(x) =  |x - v_i|$, where $x\in \R$ ($d=1$), and  $v_i = 2 \times (i \;\text{mod}\; 2) - 1$. Note that the cost functions are convex (but not strongly convex) and non-differentiable at some points. 
We exploit the used time-varying graph in~[18] where the mixing weight matrices are given by
\begin{align*}
    [W(t)]_{ij} = \begin{cases}
    \frac{r_{\langle j \rangle}}{r_{\langle t \rangle} + r_{\langle t+1 \rangle}} & i,j \in \{\langle t\rangle, \langle t+1\rangle \} \\
    1 & i = j \notin \{\langle t\rangle, \langle t+1\rangle \} \\
    0 & \text{otherwise}.
    \end{cases}
\end{align*} 
Here $\langle i \rangle = (i-1\; \text{mod}\; n) + 1$. We assume the elements of the noise vector $E(t)$ to be i.i.d.\ $\mathcal{N}(0,0.1)$ Gaussian random variables. The parameters of the dynamics in~(1) are fine-tuned to $\alpha_0 = 0.0055$ and $\beta_0 = 0.21$. To complete the validation, we implement the one time-scale without any damping mechanism for the noise vectors, i.e., $\beta(t) = 1$ for every $t$.

\textbf{Results.} In Figure~\ref{fig:sim}, top-left and top-right plots demonstrate the trajectories vs. training time for the dynamics~(1) with choices of $(\mu,\nu) = (0.6,0.77) \in \mathcal{R}_1$, and $ (\mu,\nu) =(0.2,0.3)\in \mathcal{R}_3$, respectively. It can be verified that the state of each node converges to an optimal point for the \textit{first} pair (in $\mathcal{R}_1$) of $(\mu,\nu)$ while the trajectories follow a random walk and do not converge for the \textit{second} pair (in $\mathcal{R}_3$) of $(\mu,\nu)$. In Figure~\ref{fig:sim}, the bottom-left and bottom-right plots show the trajectories vs. training time for the one time-scale method with $\nu=0.77$ and $0.3$, respectively. It can be observed that the states of the node form a random walk, and it shows the privilege of exploiting the two-time scale in the presence of noise.

\begin{figure}[h]
    \centering
    \includegraphics[width=0.45\textwidth]{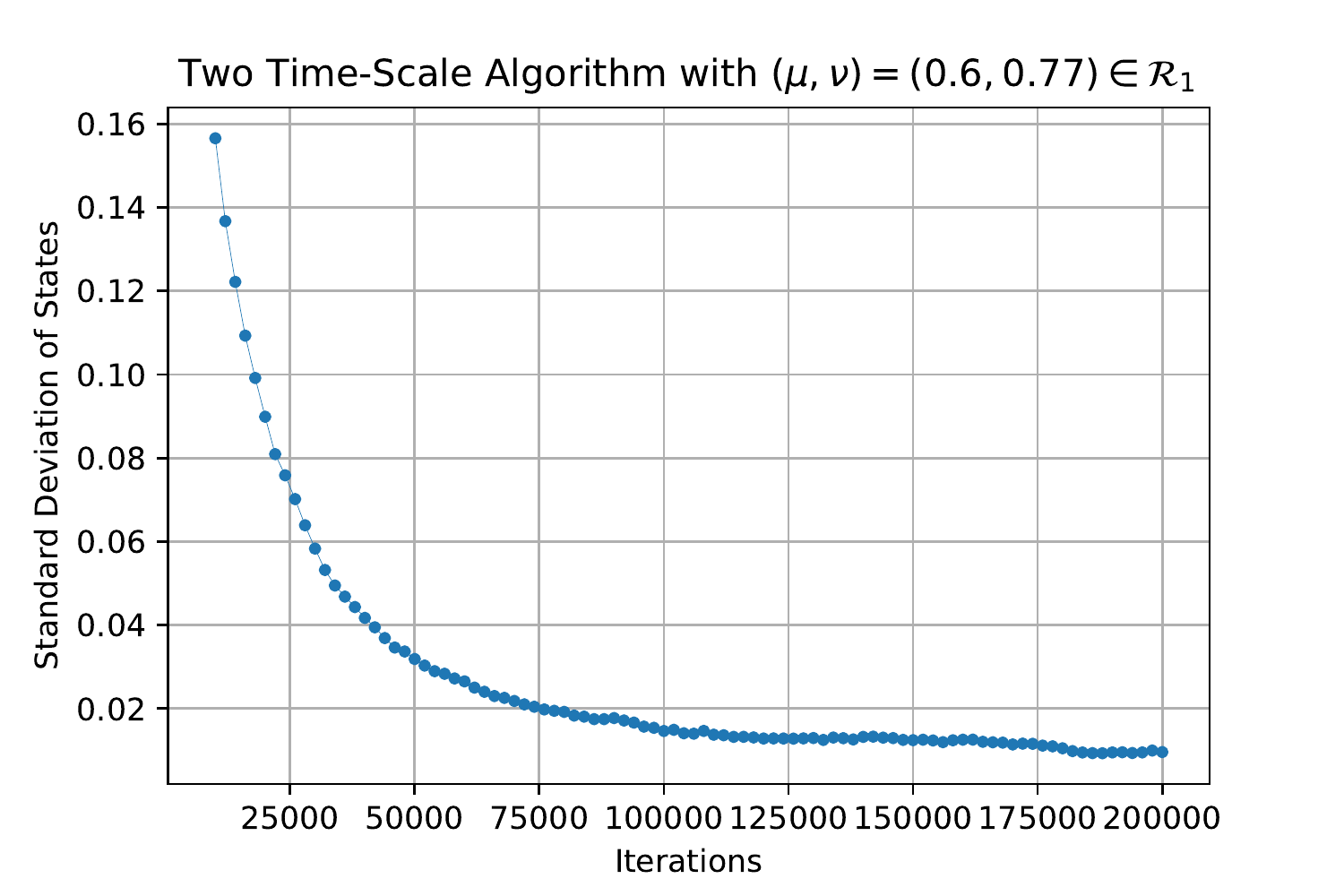}
    \vspace{-1mm}
\includegraphics[width=0.45\textwidth]{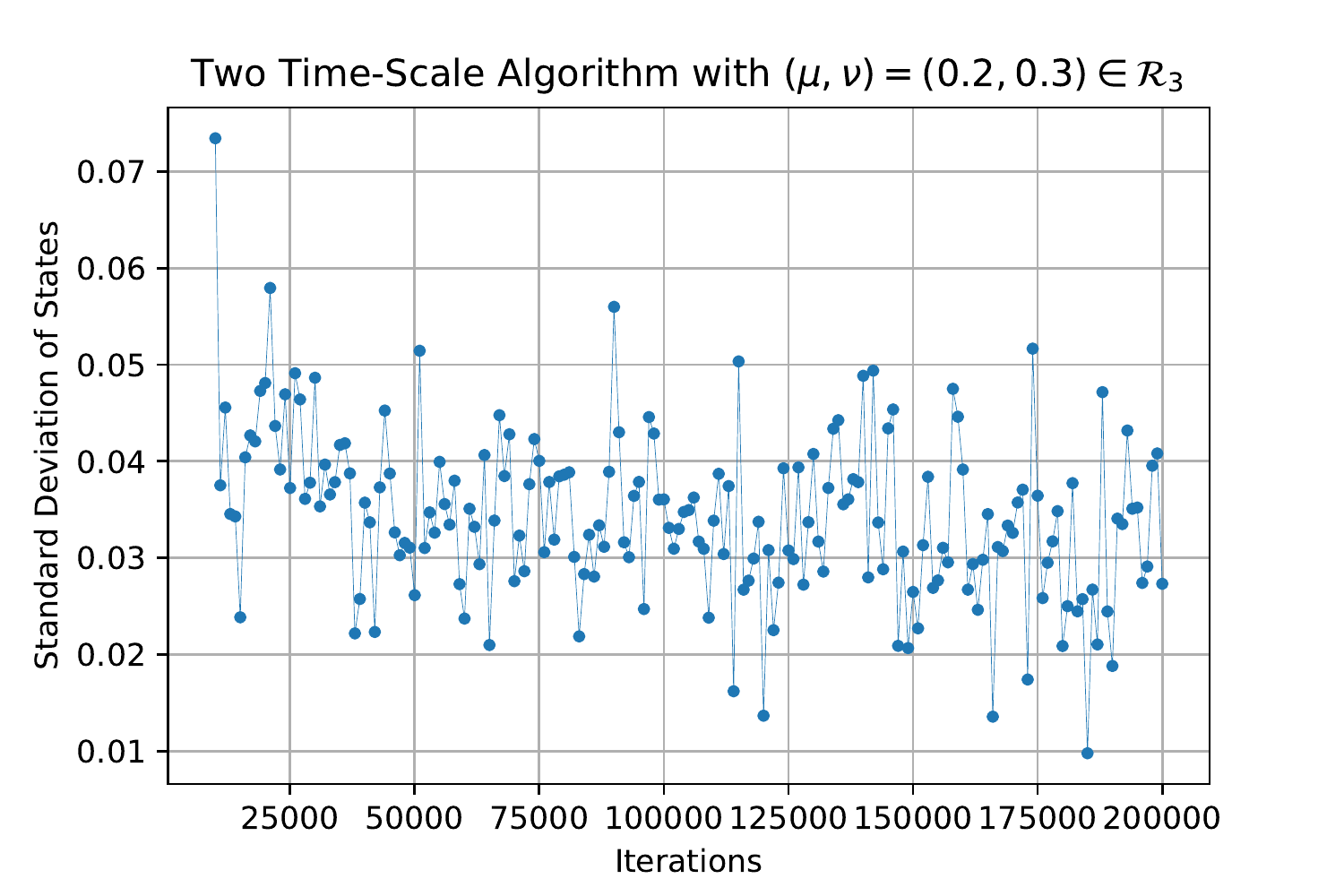}
    \includegraphics[width=0.45\textwidth]{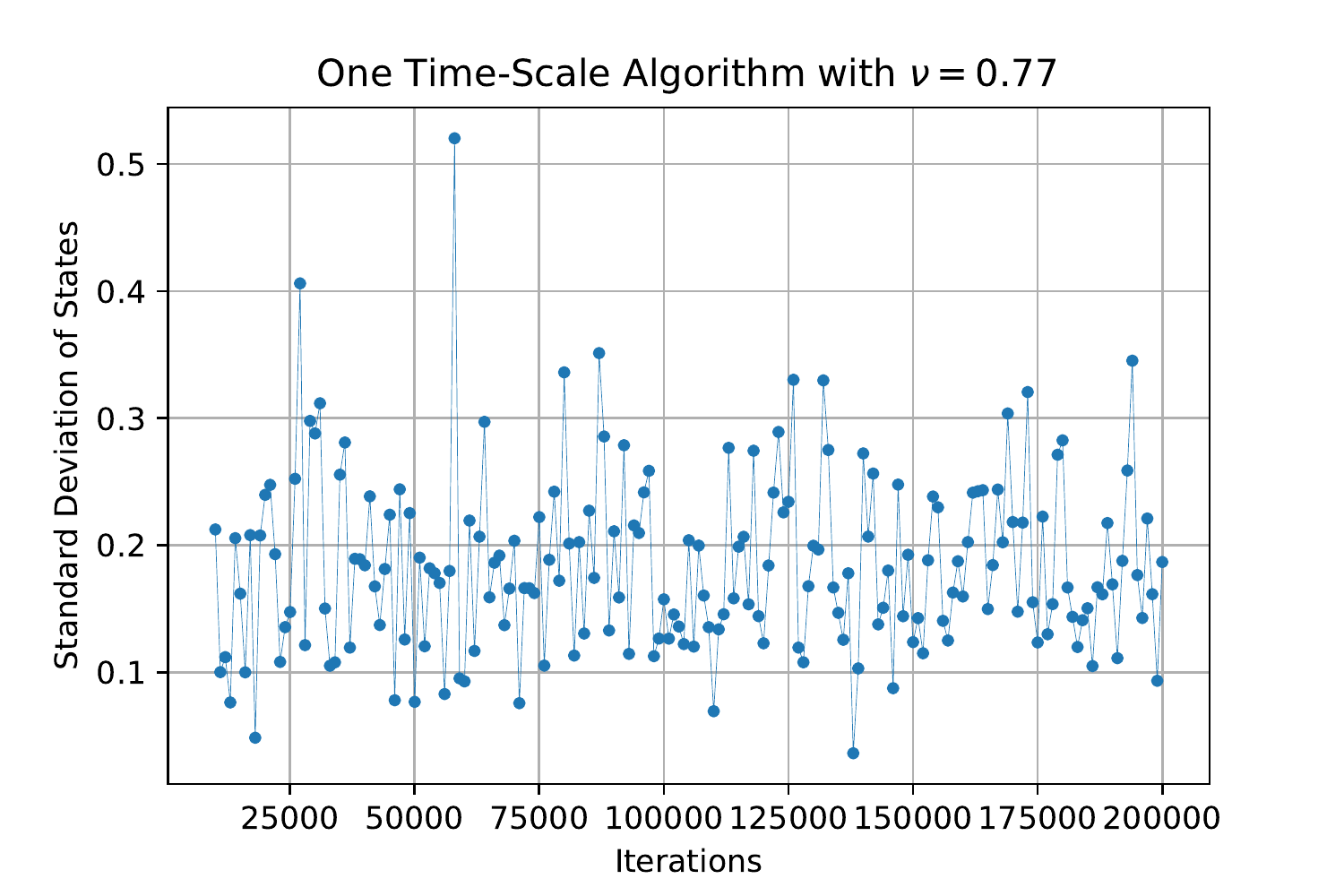}
    \vspace{-1mm}
\includegraphics[width=0.45\textwidth]{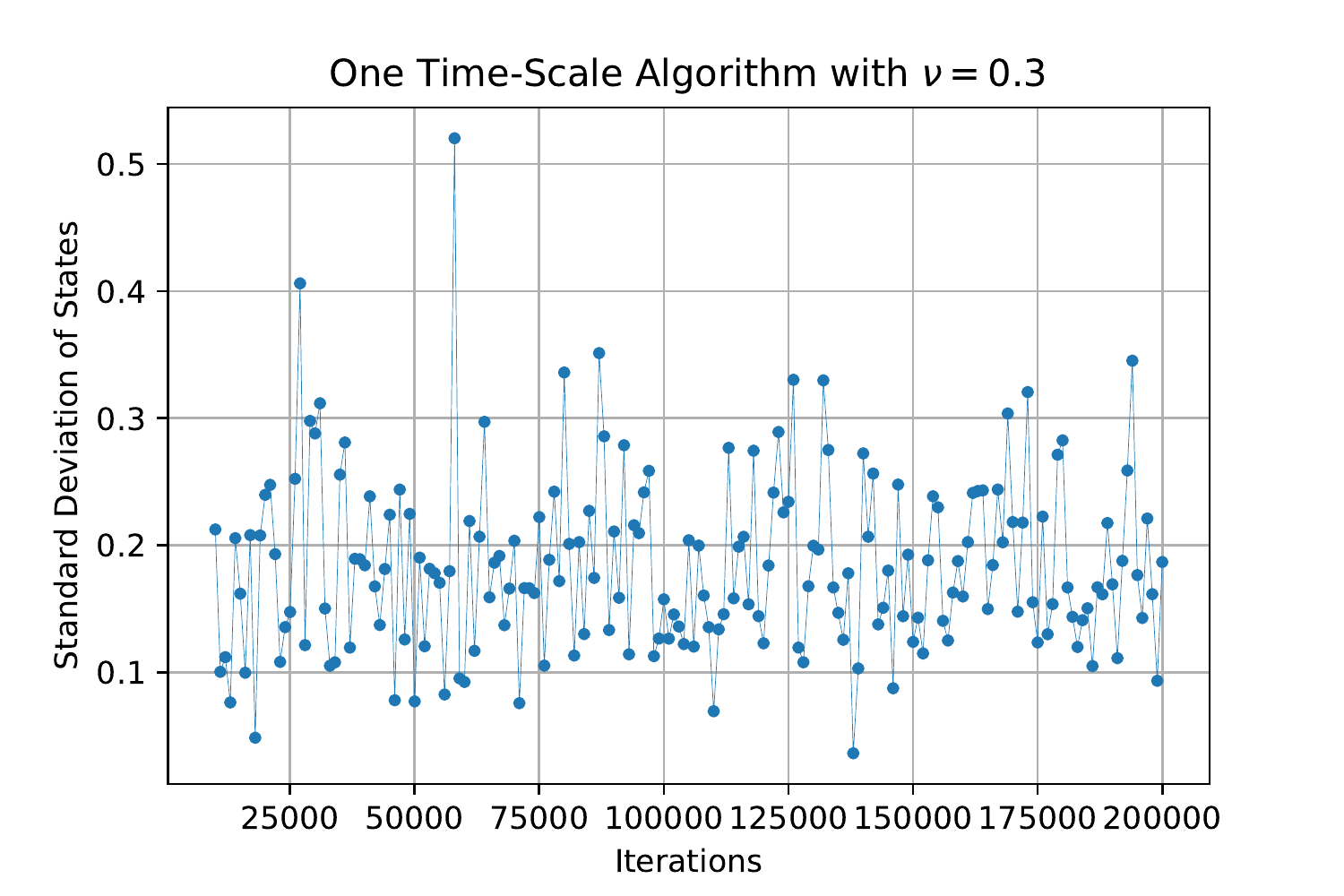}
    
    \caption{Standard Deviation of States vs. Iterations: Two Time-Scale Algorithm with $(\mu,\nu) = (0.6,0.77)$ (top-left), $(0.2,0.3)$ (top-right), and One Time-Scale Algorithm with $\nu=0.77$ (bottom-left), $0.3$ (bottom-right).}
    \label{fig:sm1}
\end{figure}

Furthermore, Figure~\ref{fig:sm1} demonstrates the deviation of nodes' states from the mean state for every iteration. It can be verified that one-time scale and two-time scale algorithms with $(\mu,\nu)$ in $\mathcal{R}_3$ both exhibit a random walk behavior for the states of the node. On the other hand, the two-time scale algorithm with $(\mu,\nu)$ in $\mathcal{R}_1$ results in consensus among the states of all nodes.

\begin{figure}
    \centering
    \includegraphics[width=0.45\textwidth]{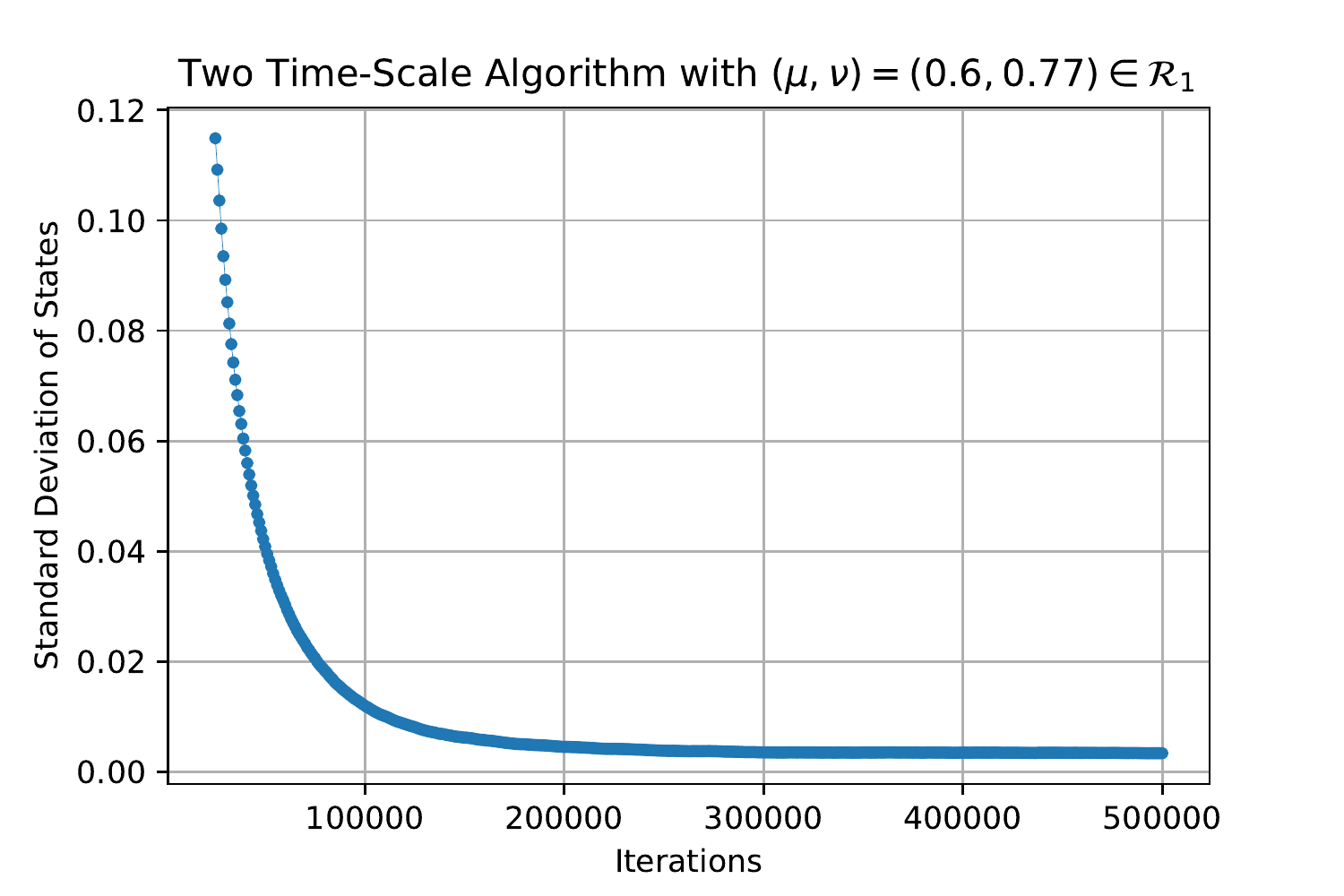}
    \includegraphics[width=0.45\textwidth]{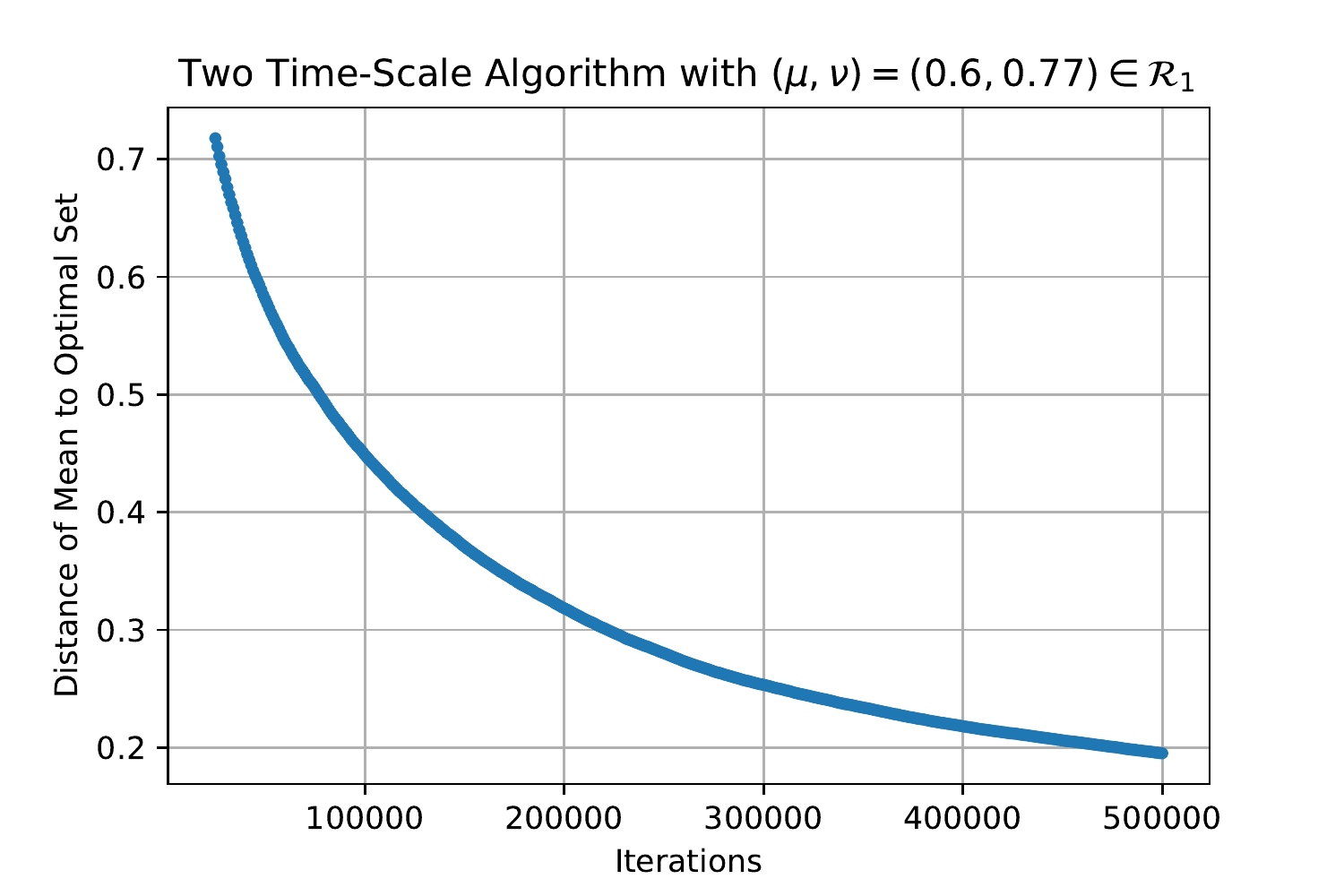}
    \caption{Standard Deviation of States and Distance of Mean State to Optimal Set vs. Iterations: Two Time-Scale Algorithm with $(\mu,\nu) = (0.6,0.77)$.}
    \label{fig:sm3}
\end{figure}

To further demonstrate the effectiveness of the two-time scale algorithm, we conducted an experiment on a time-varying network consisting of $n=6$ agents with each node's cost function defined as $f_i(\bfx) = \|\bfx - \bfv_i\|_1$, where $\bfv_i$ is a constant vector in $\mathbb{R}^{10}$. For our experiment, we set $\bfv_i = \mathbf{w}_1$ for nodes $i=1,3,5$, and $\bfv_i = \mathbf{w}_2$ for nodes $i=2,4,6$ where the elements of $\mathbf{w}_1$ and $\mathbf{w}_1$ are generated randomly from Gaussian distribution with $\mathcal{N}(0,0.1)$. It is worth noting that these cost functions are convex, but not strongly convex, and are non-differentiable at certain points. For the noise vectors, we sample points from $\mathcal{N}(0,0.1)$. The parameters of the used dynamics are fine-tuned to $\alpha_0 = 0.0075$ and $\beta_0 = 0.12$. In Figure~\ref{fig:sm3}, the left and right plots demonstrate the deviation of states' nodes from the average state and the distance of the mean state from the optimal set for every iteration. By analyzing the plots, we can observe that each node's state gradually converges to an optimal point, as the distance from the average state and the optimal set decreases over time.
  
\section{The Preliminaries}\label{sec:prelims}
In this section, we present some results which will be used in the proof of the main theorem. The proof of Lemma~\ref{lemma:transition2} and Lemma~\ref{lm:sum_pro_g} are provided in Appendix. We refer to the cited references for the proof of other preliminaries. 

The first result is known as the Robbins-Siegmund's Theorem~\cite{robbins1971convergence} as described below. 
\begin{theorem}\label{thm:robbin_siegmund}
Suppose that for non-negative random processes $\{v(t)\},\{\zeta(t)\},\{u(t)\},$ and $\{z(t)\}$ that are adapted to a filtration $\{\F_t\}$, we have
    \begin{align}\label{eq:RS-cond}
        \E[v(t+1)\mid\F_t] \leq (1 + \zeta(t))v(t) - u(t) + z(t),
    \end{align}
    almost surely, for all $t\geq 0$. Then, if $\sum_{t=1}^{\infty}z(t) < \infty$ and $\sum_{t=1}^{\infty}\zeta(t) < \infty$ almost surely, we almost surely have $\sum_{t=1}^{\infty} u(t) < \infty$ and $v(t)$ converges almost surely to a (non-negative) random variable $v$.
\end{theorem}
Inspired by~\cite[Lemma~1]{reisizadeh2021distributed}, we provide  important properties on the product of the weight matrices $\{A(t)\}_t$ in the next lemma.
\begin{lemma}\label{lemma:transition2}
    Let $\Wt$ satisfy the connectivity Assumption~\ref{Assumption:Graph} with parameters $(B,\minW)$,  and let  $\{A(t)\}$ be given by ${A(t)=(1-\beta(t))I+\beta(t)W(t)}$ where $\beta(t)\in (0,1]$ for all $t$, and $\{\beta(t)\}$ is a non-increasing sequence. Then, for any matrix $U\in \mathbb{R}^{n\times d}$, and all $s\geq 1$, we have
  \begin{align}\label{eq:lm-tran-1}
        &\Nr{\lp \Phi(t:s)-\ones \bfr^T\rp U }^2 \leq \Nr{U}^2, 
    \end{align}
    for every $t>s$. Furthermore, we have
    \begin{align}\label{eq:lm-tran-2}
        \!\!\Nr{\lp \Phi(s\!+\!B\!+\!1\!:\!s) \!-\!\ones \bfr^T\rp U }^2 \!\leq\! (1-\lambda B\beta(s\!+\!B))\Nr{U}^2\!.
    \end{align}
\end{lemma}
The proof of Lemma~\ref{lemma:transition2} is presented in Appendix.

The following theorem from~\cite{verma2022maximal} is a consequence of the Robbins-Siegmund's Theorem and  plays a crucial role in the proof of Theorem~\ref{thm:almost_sure}. 
\begin{theorem}\cite[Lemma~3]{verma2022maximal}\label{thm:a_s_conv}
Let the optimal set ${\mathcal{X}^\star=\arg \min_{\bfx\in \mathbb{R}^d} f(\bfx)}$ be nonempty for a convex and continuous function $f: \R^d \rightarrow \R$. Moreover, assume $\{\bfy(t)\}$ is a sequence satisfying 
\begin{align*}
    \EE{\|\bfy(t\!+\!1)-\bfx^\star \|^2 | \F_t} \leq & (1\!+\!\zeta(t)) \|\bfy(t)-\bfx^\star \|^2 \\
    &\phantom{\leq} - \xi(t) (f(\bfy(t))\!-\!f(\bfx^\star))\!+\!z(t),
\end{align*}
for all $t\geq 1$ and for all $\bfx^\star\in \mathcal{X}^\star$ almost surely, where non-negative sequences  ${\zeta(t)}$, ${\xi(t)}$, and ${z(t)}$ for satisfy $\sum_{t=1}^{\infty}\zeta(t)< \infty$, $\sum_{t=1}^{\infty}\xi(t)=\infty$, and $\sum_{t=1}^{\infty}z(t)<\infty$. Then, the sequence $\{\bfy(t)\}$ converges to some solution $\tilde{\bfx} \in \mathcal{X}^\star$ almost surely.
\end{theorem}

To show the consensus in the almost sure sense, we exploit the following result which has been proven as part of~\cite[Theorem~1]{reisizadeh2021distributed}.
\begin{lemma}\label{lm:consen2}
Let Assumptions~\ref{Assumption:Noise},~\ref{Assumption:Graph}, and~\ref{Assumption:Function}-(c) are satisfied for the dynamics~\eqref{eqn:UpdateLawMatrixForm}. Then for  non-increasing sequence $\{\beta(t)\}$ with $0<\beta(t)\leq 1$ for all $t\geq 1$, we have
\begin{align*}
	   \EE{\Nr{X(t)\!-\!\ones\bar{\bfx}(t)}^2} &\leq  \cThree\sum_{s=1}^{t-1} \lb \beta^2(s)\! \prod_{k=s+1}^{t-1} (1\!-\!\lambda \beta(k)) \rb\\
	   & \phantom{\leq} +\! \cFour \sum_{s=1}^{t-1}\lb\frac{\alpha^2(s)}{\beta(s)}\! \prod_{k=s+1}^{t-1}\!\!\lp 1\!-\!\lambda \beta(k)\rp^{\frac{1}{2}}\rb\!,
	\end{align*}
\end{lemma}
for some constants $\cThree,\cFour>0$.
\begin{lemma}\label{lm:sum_pro_g}
        Let $\{p(t)\}$ and $\{q(t)\}$ be two positive and non-increasing sequences for $t\geq 1$ and there exists for $0<A<1$ such that 
        \begin{align}\label{eq:sum_prod_con}
           -\Delta p(t) \leq A p(t) q(t), 
        \end{align}
        for every $t\geq t_0$.
        Then, we have
        \begin{align}\label{eq:sum_dl_m}
            \sum_{s=1}^{t-1}  \lb p(s) \prod_{k=s+1}^{t-1} \lp 1-q(k) \rp \rb &\leq S\frac{p(t)}{q(t)}.
        \end{align}
for every $t\geq t_0$, and some positive $S$ which is not a function of $t$ (but may depend on $\tz$ and $A$). 
\end{lemma}
\begin{lemma}\cite[Lemma~3]{reisizadeh2021distributed}\label{lem:CauchyExtn}
For any pair of vectors $u,v$ and any scalar $\theta > 0 $, we have 
\begin{align*}
    \norm{u+v}^2 \leq (1+\theta)\norm{u}^2 + \lp 1+\frac{1}{\theta}\rp\norm{v}^2.
\end{align*}
Similarly, for matrices, $U$ and $V$, and a scalar $\theta > 0 $, we have
\begin{align*}
    \norm{U+V}^2_\bfr \leq (1+\theta)\norm{U}^2_\bfr + \left(1+\frac{1}{\theta}\right)\norm{V}^2_\bfr.
\end{align*}

\end{lemma}
\begin{lemma}\cite[Lemma~3]{reisizadeh2022dimix}\label{lm:sum-exp}
    For any $\dl\in \mathbb{R}$, $\tau \geq 0$, and $T\geq 1$, we have
    \begin{align}
    \sum_{t=1}^T (t+\tau)^{\dl} \leq
    \begin{cases}
    \frac{\tau^{1+\delta}}{|1+\dl|} & \textrm{if $\dl<-1$},\\
    \ln\lp \frac{T}{\tau} +1 \rp& \textrm{if $\dl=-1$},\\
    \frac{2^{1+\dl}}{1+\dl} (T+\tau)^{1+\dl} & \textrm{if $\dl>-1 $}.
    \end{cases}
    \end{align}
    \end{lemma}
\section{Proof of Theorem~\ref{thm:almost_sure}}\label{proof:main_result}
In this section, we provide the proof of Theorem~\ref{thm:almost_sure}. We first show that the deviation of the agents' states from their \textit{average} converges to a random variable, which is later shown to be zero with probability $1$.  Next, we analyze the distance of the average state from an arbitrary point $\bfx^\star$  in the optimal set $\mathcal{X}^\star$ of function $f(\cdot)$. These together lead to the proof of the theorem.

\subsection{State Deviation from the Average State}
\label{sec:proof:variance}
We first define $\dl(t)\!:=\!\Nr{\Dl(t)}$ where ${\Dl(t)\!:=\!X(t)\!-\!\ones\mx(t)}$ and $\mx(t) :=\bfr^T X(t) = \sum_{i=1}^n r_i \bfx_i(t)$. Our ultimate goal is to show that $\dl(t)$ vanishes almost surely. To that end, we first show its convergence in this section, and then show that it converges to $0$ in Section~\ref{sec:state-to-average}. Since we are dealing with time-varying graphs, we cannot guarantee any decent for $\dl(t)$ in every iteration. However, since the graph is {$B$-connected} (see Assumption~\ref{Assumption:Graph}), such a claim can be made from $\dl(t)$ and $\dl(t+B)$. As a result, we can show that for any $1\leq \tau \leq B$ the sequence $\lc\dl(\tau \!+\! k B)\rc_{k=0}^{\infty}$ converges almost surely to a (non-negative) random variable $v_\tau$.

Starting from~\eqref{eqn:UpdateLawMatrixForm}, we can write
\begin{align}\label{eq:x-t-mat}
	    X(t) = \sum_{s=1}^{t-1}\Phi(t:s)U(s) + \Phi(t:0)X(1).
\end{align}
Assuming $X(1)={\mathbf{0}}$, the dynamics in~\eqref{eq:x-t-mat} reduces to 
\begin{align}\label{eq:x-t-mat-2}
	    X(t) = \sum_{s=1}^{t-1}\Phi(t:s)U(s).
\end{align}
By multiplying both sides of~\eqref{eq:x-t-mat-2} from the left by $\bfr^T$ and using the fact $\bfr^T A(t)=\bfr^T$ for every $t\geq 1$, we get $\mx = \bfr^T X(t) = \sum_{s=1}^{t-1} \bfr^T \Phi(t:s) U(s) = \sum_{s=1}^{t-1} \bfr^T U(s)$. Subtracting $\ones \mx(t)$ from~\eqref{eq:x-t-mat-2}, we get
\begin{align*}
    \Dl(t) &= \sum_{s=1}^{t-1}(\Phi(t:s) - \ones\bfr^T)U(s).
\end{align*}
Writing this equation for iteration $\k+B$ we get
\begin{align}\label{eq:dl-t-1}
\Dl(\k+B) & =  \sum_{s=1}^{\k+B-1}(\Phi(\k+B:s) - \ones\bfr^T)U(s)\nonumber\\
    & = \sum_{s=1}^{\k-1}(\Phi(\k+B:s) - \ones\bfr^T)U(s) + \sum_{s=\k}^{\k+B-1}(\Phi(\k+B:s) -\ones\bfr^T)U(s)\nonumber\\
    & \stackrel{\rm{(a)}}{=} \Phi(\k+B:\k\!-\!1)\sum_{s=1}^{\k-1}(\Phi(\k:s) - \ones\bfr^T)U(s) + \sum_{s=\k}^{\k+B-1}(\Phi(\k+B:s) -\ones\bfr^T)U(s)\nonumber\\
    & = \!\Phi(\k\!+\!B:\k\!-\!1)\Dl(\k) \!+\! \sum_{s=\k}^{\k+B-1}(\Phi(\k\!+\!B\!:\!s) \!-\!\ones\bfr^T)U(s),
\end{align}
where in~$\rm{(a)}$ we used the fact $A(s)\ones=\ones$ for every $s\geq 1$. This leads to
\begin{align}\label{eq:dl-t}
    \EE{ \dl^2(\k\!+\!B)\md\F_{\k}} 
    &= \EE{ \Nr{D(\k\!+\!B)}^2\md\F_{\k}} \nonumber\\
    & = \mathbb{E}\Bigg[\Biggl|\!\Biggl|\Phi(\k+B:\k\!-\!1)\Dl(\k) +\! \sum_{s=\k}^{\k+B-1} (\Phi(\k\!+\!B\!:\!s)\!-\! \ones\bfr^T)U(s)\Biggr|\!\Biggr|_{\bfr}^2 \Bigg\vert \F_{\k}\Bigg]\nonumber \\
    &= \EE{\Nr{\Phi(\k+B:\k-1)\Dl(\k)}^2 \md \F_{\k}} + \EE{\Nr{\sum_{s=\k}^{\k+B-1} (\Phi(\k+B:s)- \ones\bfr^T)U(s)}^2 \md \F_{\k}} \nonumber \\ 
    &\phantom{=} + 2 \sum_{i=1}^{n}r_i \mathbb{E}\Bigg[ \Bigg<\sum_{s=\k}^{\k+B-1} \lb (\Phi(\k+B:s)- \ones\bfr^T)U(s)\rb_i ,\lb \Phi(\k+B:\k-1)\Dl(\k)\rb_i\Bigg > \Bigg\vert \F_{\k} \Bigg].
\end{align}
Next, we bound each term in~\eqref{eq:dl-t}, separately. Note that $\bfr^T X(t) = \mx(t)$ and $\bfr^T \ones = 1$. Hence 
\begin{align*}
    \ones \bfr^T \Dl(t) =  \ones \bfr^T (X(t) - \ones \mx(t)) = \ones \mx(t) - \ones \mx(t) =0.
\end{align*}
Therefore, using Lemma~\ref{lemma:transition2},  the first term can be bounded  as
\begin{align}\label{eq:T1}
   \EE{\Nr{\Phi(\k+B:\k-1)\Dl(\k)}^2 \md \F_{\k}}
   & = \EE{\Nr{\lp \Phi(\k+B:\k-1) - \ones \bfr^T \rp \Dl(\k)}^2 \md \F_{\k}}\nonumber\\
   & = \Nr{\lp \Phi(\k+B:\k-1) -\ones \bfr^T\rp\Dl(\k)}^2\nonumber\\
   & \leq  \lp 1- \lambda B \beta(\k+B-1) \rp \dl^2(\k).
\end{align}
Next, in order to bound the second term, we can use the convexity of the $\Nr{\cdot}$ to write
\begin{align}\label{eq:T2-1}
&\EE{\Nr{ \sum_{s=\k}^{\k+B-1} (\Phi(k+B:s)- \ones\bfr^T)U(s)}^2 \md \F_{\k} } \nonumber\\
& \leq  \EE{ B \sum_{s=\k}^{\k+B-1}  \Nr{(\Phi(\k+B:s)- \ones\bfr^T)U(s)}^2 \md \F_{\k} }\nonumber \\
&\leq 2B\sum_{s=\k}^{\k+B-1} \beta^2(s)\EE {\Nr{ (\Phi(\k+B:s) \!-\! \ones\bfr^T)E(s)}^2 \md \F_{\k} } \nonumber \\
&\phantom{\leq} \!+\! 2B\!\sum_{s=\k}^{\k+B-1}\!\!\alpha^2(s)\EE{\Nr{ (\Phi(\k\!+\!B\!:\!s) \!-\! \ones\bfr^T)\G(s)}^2 \md \F_{\k} }\!,
\end{align}
where the second inequality follows from the fact that $U(s) = \beta(s) E(s) - \alpha(s) \G(s)$ and  Lemma~\ref{lem:CauchyExtn} with $\theta=1$.
For the first term in~\eqref{eq:T2-1}, from Lemma~\ref{lemma:transition2} we have
\begin{align}\label{eq:T2-2}
    &\sum_{s=\k}^{\k+B-1}\beta^2(s)\EE{\Nr{ (\Phi(\k+B:s) \!-\! \ones\bfr^T)E(s))}^2 \md \F_{\k} } \nonumber\\
    & \leq \sum_{s=\k}^{\k+B-1}\beta^2(s)\EE{\Nr{E(s)}^2\md \F_{\k}} \leq \gamma \sum_{s=\k}^{\k+B-1} \beta^2(s),
\end{align}
where the last inequality follows from Assumption~\ref{Assumption:Noise} for every $s\geq \k$ we can write
\begin{align*}
    \EE{\Nr{E(s)}^2\md \F_{\k}} & = \EE{\EE{\Nr{E(s)}^2\md \F_{s}}\md \F_{\k}}\nonumber\\
    & = \sum_{i=1}^{n}r_i \EE{\EE{\|\bfe_i(s)\|^2\md \F_{s}}\md \F_{\k}} \nonumber\\
    & \leq \sum_{i=1}^{n}r_i\EE{\gamma \md \F_{\k}} = \gamma.
\end{align*}
Similarly, for the second term in~\eqref{eq:T2-1}, we can apply Lemma~\ref{lemma:transition2} and write
\begin{align}\label{eq:T2-3}
     \sum_{s=\k}^{\k+B-1} \alpha^2(s)\EE{\Nr{ (\Phi(\k+B:s) \!-\! \ones\bfr^T)\G(s)}^2 \md \F_{\k}} & \leq \sum_{s=\k}^{\k+B-1} \alpha^2(s) \EE{\Nr{\G(s)}^2\md \F_{\k}} \nonumber\\
     & \leq L^2 \sum_{s=\k}^{\k+B-1} \alpha^2(s),
\end{align}
where in the last inequality follows from Assumption~\ref{Assumption:Function}-(c).
Plugging~\eqref{eq:T2-2} and~\eqref{eq:T2-3} into~\eqref{eq:T2-1}, we arrive at
\begin{align}\label{eq:T2}
    \EE{\Nr{ \sum_{s=\k}^{\k+B-1} (\Phi(\k+B:s)- \ones\bfr^T)U(s)}^2 \md \F_{\k} }
     & \leq 2B \sum_{s=\k}^{\k+B-1} \lp \gamma\alpha^2(s) + L^2 \beta^2(s)\rp.
\end{align}
Next, we focus on each term of the last summation  in~\eqref{eq:dl-t}. Since $U(s) = \beta(s) E(s) - \alpha(s) \G(s)$, we have
\begin{align}\label{eq:T3-0}
    \mathbb{E} & \Bigg[ \Bigg<\sum_{s=\k}^{\k+B-1} \lb (\Phi(\k+B:s)- \ones\bfr^T)U(s)\rb_i, \lb \Phi(\k+B:\k-1)\Dl(\k)\rb_i\Bigg > \Bigg\vert \F_{\k} \Bigg]\nonumber\\
    & = \mathbb{E}\Bigg[ \Bigg<\sum_{s=\k}^{\k+B-1} \lb (\Phi(\k+B:s)- \ones\bfr^T)\beta(s)E(s)\rb_i, \lb \Phi(\k+B:\k-1)\Dl(\k)\rb_i\Bigg > \Bigg\vert \F_{\k} \Bigg]\nonumber\\
    & \phantom{=} + \!\mathbb{E}\Bigg[\! \Bigg<\!\sum_{s=\k}^{\k+B-1} \!\!\lb  (\Phi(\k\!+\!B\!:\!s) \!-\! \ones\bfr^T)(-\alpha(s)\G(s))\rb_i\!, \lb \Phi(\k+B:\k-1)\Dl(\k)\rb_i\Bigg > \Bigg\vert \F_{\k} \Bigg].
\end{align}
For the first term in~\eqref{eq:T3-0}, we have
\begin{align}\label{eq:T3-1}
    \mathbb{E} & \Bigg[ \Bigg<\sum_{s=\k}^{\k+B-1} \lb (\Phi(\k+B:s)- \ones\bfr^T)\beta(s)E(s)\rb_i, \lb \Phi(\k+B:\k-1)\Dl(\k)\rb_i\Bigg > \Bigg\vert \F_{\k} \Bigg]\nonumber\\
    & = \mathbb{E} \Bigg[ \!\Bigg<\sum_{s=\k}^{\k+B-1} \beta(s) \lb (\Phi(\k\!+\!B\!:\!s)- \ones\bfr^T)\EE{ E(s)\md \F_{\k}}\rb_i, \lb  \Phi(\k\!+\!B\!:\!\k-1)\Dl(\k)\rb_i\Bigg > \Bigg] = 0,
\end{align}
where the last inequality holds since from Assumption~\ref{Assumption:Noise} for every $s\geq \k$ we get
\begin{align*}
    \EE{ E(s)\md \F_{\k}} = \EE{ \EE{E(s)\md \F_s}\md \F_{\k}} = 0.
\end{align*}
For the second term in~\eqref{eq:T3-0}, using the Cauchy-Schwarz inequality and the fact that $2ab\leq a^2\!+\!b^2$, we can write
    \begin{align}\label{eq:T3-2}
    \mathbb{E}&\Bigg[ \Bigg<\sum_{s=\k}^{\k+B-1} \!\lb (\Phi(\k\!+\!B\!:\!s) \!-\! \ones\bfr^T)(-\alpha(s)\G(s))\rb_i, \lb \Phi(\k+B:\k-1)\Dl(\k)\rb_i\Bigg > \Bigg\vert \F_{\k} \Bigg]\nonumber\\
    & \leq \mathbb{E} \Bigg[\Biggl|\!\Biggl|  \sum_{s=\k}^{\k+B-1} \!\!\alpha(s) \lb (\Phi(\k\!+\!B\!:\!s) \!-\! \ones\bfr^T)\G(s)\rb_i \Biggr|\!\Biggr| \times \lnr \lb \Phi(\k+B:\k-1)\Dl(k)\rb_i \rnr \vert \F_{\k} \Bigg]\nonumber\\
    & = \frac{1}{2} \mathbb{E} \Bigg[ \frac{1}{\omega(\k)} \Biggl|\!\Biggl| \sum_{s=\k}^{\k+B-1}\!\!\alpha(s)\lb (\Phi(\k\!+\!B\!:\!s)\!-\! \ones\bfr^T) \G(s)\rb_i \Biggr|\!\Biggr|^2 + \omega(\k)\lnr \lb \Phi(\k+B:\k\!-\!1)\Dl(\k)\rb_i \rnr^2 \vert \F_{\k} \Bigg],
\end{align}
for any $\omega(\k)>0$, which will be determined later. Hence, using~\eqref{eq:T3-1} and~\eqref{eq:T3-2} in~\eqref{eq:T3-0} we get
\begin{align}\label{eq:T3}
    & 2 \sum_{i=1}^{n}r_i \mathbb{E}\Bigg[ \Bigg<\sum_{s=\k}^{\k+B-1} \lb (\Phi(\k+B:s)- \ones\bfr^T)U(s)\rb_i, \lb \Phi(\k+B:\k-1)\Dl(\k)\rb_i\Bigg > \Bigg\vert \F_{\k} \Bigg] \nonumber\\
    &  \leq  \frac{1}{\omega(\k)} \sum_{i=1}^{n}\hspace{-1pt} r_i \mathbb{E}\left[  \left\| \sum_{s=\k}^{\k+B-1} \alpha(s)  \right.\left[ (\Phi(\k + B : s)  -  \ones\bfr^T) \G(s)\hspace{-1pt} \right]_i \right\|^2  \Big| \F_{\k}\Bigg]\nonumber\\
    & \phantom{\leq} +  \omega(\k) \sum_{i=1}^{n} r_i \lnr \lb \Phi(\k+B:\k-1)\Dl(\k)\rb_i \rnr^2\nonumber \\
    & \!\stackrel{\rm{(a)}}{\leq} \!\!\frac{B}{ \omega(\k)}\!\!\sum_{s=\k}^{\k+B-1} \!\!\alpha^2(s)\mathbb{E}\hspace{-1pt}\left[ \!\Nr{(\Phi(\k\!+\!B\!:\!s) \!-\! \ones\bfr^T)\G(s)}^2 \hspace{-1pt}\md \F_{\k} \hspace{-1pt}\right] +  \omega(\k)  \Nr{\Phi(\k+B:\k\!-\!1)\Dl(\k)}^2\nonumber\\
    & \!\stackrel{\rm{(b)}}{\leq}\! \frac{B L^2}{ \omega(\k)}  \sum_{s=\k}^{\k+B-1}\!\alpha^2(s) \!+\!   \omega(\k)\!  \lp 1 \!-\! \lambda B \beta(\k\!+\!B\!-\!1) \rp\dl^2(\k),
\end{align}
where  step $\rm{(a)}$ follows from the convexity of $\Nr{\cdot}$ and the inequality in $\rm{(b)}$ follows from~\eqref{eq:T2-3} and~\eqref{eq:T1}.
Finally, plugging~\eqref{eq:T1},~\eqref{eq:T2}, and~\eqref{eq:T3} into~\eqref{eq:dl-t} we have
\begin{align}\label{eq:robin-dl}
    \EE{\dl^2(\k+B)\md\F_{\k}} & \leq \lp 1+ \omega(\k) \rp\lp 1- \lambda B \beta(\k+B-1) \rp \dl^2(\k) \nonumber\\
    &\phantom{=} + 2B \sum_{s=\k}^{\k+B-1} \lp \gamma\alpha^2(s) + L^2 \beta^2(s)\rp + \frac{B L^2}{\omega(\k)} \sum_{s=\k}^{\k+B-1}\alpha^2(s)\nonumber\\
    & \stackrel{\rm{(a)}}{\leq} \lp 1+ \omega(\k)  \rp\lp 1- \lambda B \beta(\k+B-1) \rp \dl^2(\k) \nonumber\\
    &\phantom{=} + 2B^2 \lp \gamma\alpha^2(\k) + L^2 \beta^2(\k)\rp + \frac{B^2 L^2}{ \omega(\k)} \alpha^2(\k) \nonumber\\
    & \stackrel{\rm{(b)}}{=} \lp 1+\lambda B\beta(\k) \rp\lp 1- \lambda B \beta(\k+B-1) \rp \dl^2(\k) \nonumber\\
    &\phantom{=} + 2B^2 \lp \gamma\alpha^2(\k) + L^2 \beta^2(\k)\rp + \frac{B L^2}{\lambda} \frac{\alpha^2(\k)}{\beta(\k)}\nonumber\\
    & \leq \lp 1+\lambda B (\beta(\k) - \beta(\k+B-1)) \rp \dl^2(\k) \nonumber\\
    &\phantom{=}  + 2B^2 \lp \gamma\alpha^2(\k) + L^2 \beta^2(\k)\rp + \frac{B L^2}{\lambda} \frac{\alpha^2(\k)}{\beta(\k)},
\end{align}
where the inequality in~$\rm{(a)}$ follows from this assumption that $\{\alpha(t)\}$ and $\{\beta(t)\}$ are non-increasing step-size sequences and in the step~$\rm{(b)}$ we set $\omega(\k)= \lambda B \beta(\k)$. 

Now, for $\tau=1,2,\ldots,B$, consider $B$ random processes $\left\{\dl^2(\tau + kB)\right\}_{k=0}^\infty$. Due to~\eqref{eq:robin-dl}, each of these random processes satisfy the inequality~\eqref{eq:RS-cond} in Theorem~\ref{thm:robbin_siegmund}, with 
\begin{align*}
    &\zeta(k) := \lambda B (\beta(\tau+ kB) - \beta(\tau+  (k+1)B-1)),\\
    & u(k) := 0,\\
    & z(k) := 2B^2 \lp \gamma\alpha^2(\tau+ kB) + L^2 \beta^2(\tau+ kB)\rp + \frac{B  L^2}{\lambda} \frac{\alpha^2(\tau+ kB)}{\beta(\tau+ kB)}.
\end{align*}
Since $\{\beta(t)\}$ is a non-increasing sequence, we have
\begin{align*}
    \sum_{k=0}^{\infty} \zeta(k) & =  B\lambda\sum_{k=0}^{\infty} \beta(\tau+kB) - \beta(\tau+(k+1)B)-1) \\
    & \leq B\lambda  \sum_{k=0}^{\infty} \beta(\tau+kB)  -\beta(\tau+(k+1)B)\\ & \leq B\lambda \beta(\tau) <\infty.
\end{align*}
Moreover, we have $\sum_{k=0}^\infty u(k)=0$ and Assumption~\ref{asm:step-size}-\eqref{asm:4-a2-b2} and~\eqref{asm:4-a2/b} imply that $\sum_{k=0}^\infty z(k)<\infty$. Thus, all the
conditions of Robbin-Sigmund Theorem are satisfies, and hence, for any  $\tau=1,\ldots,B$, the random process $\left\{\dl^2(\tau + kB)\right\}_{k=0}^\infty$ converges, almost surely. Consequently, there exist random variables $v_\tau$ such that $\lim_{k\to\infty}\delta^2(\tau+kB) = v_\tau$ almost surely, for $\tau=1,\ldots,B$. 

\subsection{Accumulative Variance from the States}
\label{sec:var}

In this section, we study the summation $\sum_{t=1}^\infty \alpha(t)\dl(t)$, and show that it converges. This will imply that $\dl(t)$ converges to zero. Starting from Lemma~\ref{lm:consen2} and the fact that 
\[{\sqrt{a+b} \leq \sqrt{a} + \sqrt{b}},\] we get 
\begin{align}\label{eq:exp-dl}
     \EE{\dl(t)} & = \EE{\|X(t) - \ones\mx(t)\|_\bfr} \nonumber\\
     & \leq  \lp \cThree \sum_{s=1}^{t-1} \lb \beta^2(s) \prod_{k=s+1}^{t-1} (1-\lambda \beta(k)) \rb \rp^{\frac{1}{2}} 
     \\\nonumber &\phantom{\leq} 
     +   \lp \cFour\sum_{s=1}^{t-1}\lb\frac{\alpha^2(s)}{\beta(s)} \prod_{k=s+1}^{t-1}\lp 1-\lambda \beta(k)\rp^{\frac{1}{2}}\rb\rp^{\frac{1}{2}}\!.
\end{align}
Next, we bound each summation in~\eqref{eq:exp-dl}. To this end, we  use  Lemma~\ref{lm:sum_pro_g} with $p_1(t)=\beta^2(t)$, $q_1(t)=\lambda \beta(t)$, and $A_1 = \frac{2\cOne}{\lambda}<1$ for the first summation. Using the fact that $\{\beta(t)\}$ is a non-increasing sequence and Assumption~\ref{asm:step-size}-\eqref{asm:4-Db}, for every $t\geq \tz$ we have
\begin{align*}
    -\Delta p_1(t) & =  \beta^2(t)-\beta^2(t+1)\\
    & = -\Delta \beta(t)(\beta(t)+\beta(t+1)) \nonumber\\
    & \leq (\cOne \beta^2(t))\cdot (2\beta(t)) 
    = 2\cOne \beta^3(t) = A_1 p_1(t) q_1(t).
\end{align*}
Thus, Lemma~\ref{lm:sum_pro_g} leads  to
\begin{align}\label{eq:exp-dl:T1}
    \sum_{s=1}^{t-1} \lb \beta^2(s) \prod_{k=s+1}^{t-1} (1-\lambda \beta(k)) \rb \leq  \frac{\cfone}{\lambda} \beta(t),
\end{align}
for some constant $\cfone$ and all $t\geq \tz$. 

Similarly, we use Lemma~\ref{lm:sum_pro_g} with $p_2(t)=\frac{\alpha^2(t)}{\beta(t)}$, ${q_2(t)=\frac{\lambda}{2} \beta(t)}$, and ${A_2 = \frac{4\cTwo}{\lambda}<1}$ to bound the second summation in~\eqref{eq:exp-dl}. Using the fact that  $\{\alpha(t)\}$ and $\{\beta(t)\}$  are non-increasing sequences and Assumption~\ref{asm:step-size}-\eqref{asm:4-Da}, we can write
\begin{align*}
    - \Delta p_2(t) & =  \frac{\alpha^2(t)}{\beta(t)} - \frac{\alpha^2(t+1)}{\beta(t+1)} \\
    & \leq  \frac{\alpha^2(t)}{\beta(t)} - \frac{\alpha^2(t+1)}{\beta(t)} \\
    & = \frac{(-\Delta \alpha(t))(\alpha(t+1)+\alpha(t))}{\beta(t)} \nonumber\\
    & \leq  \frac{(\cTwo  \alpha(t) \beta(t))(2\alpha(t))}{\beta(t)} \nonumber\\
    & \leq 2\cTwo \alpha^2(t) = A_2 p_2(t) q_2(t),
\end{align*}
for $t\geq \tz$.
Thus,  Lemma~\ref{lm:sum_pro_g} together with the fact that ${\sqrt{1-x}\leq 1\!-\!x/2}$ imply
\begin{align}\label{eq:exp-dl:T2}
    \sum_{s=1}^{t-1}&\lb\frac{\alpha^2(s)}{\beta(s)}\! \prod_{k=s+1}^{t-1}\lp 1\!-\!\lambda \beta(k) \rp^{\frac{1}{2}}\rb\nonumber\\
    & \leq \sum_{s=1}^{t-1}\lb\frac{\alpha^2(s)}{\beta(s)}\! \prod_{k=s+1}^{t-1}\lp 1\!-\!\frac{\lambda}{2} \beta(k) \rp \rb \cr 
    &\leq \frac{2\cftwo}{\lambda} \frac{\alpha^2(t)}{\beta^2(t)},
\end{align}
for some constant  $\cftwo$ and every $t\geq \tz$. Plugging~\eqref{eq:exp-dl:T1} and~\eqref{eq:exp-dl:T2} into~\eqref{eq:exp-dl}, we get \begin{align} 
      \EE{\dl(t)}
     & \leq  \lp  \frac{\cThree \cfone}{\lambda} \beta(t) \rp^{\frac{1}{2}} 
     +   \lp \frac{2 \cFour \cftwo}{\lambda} \frac{\alpha^2(t)}{\beta^2(t)}\rp^{\frac{1}{2}}\!.
\end{align}
Therefore, from Equation~\eqref{eq:sum_a_b_hf} in Remark~\ref{rem:asm:step-size} and Assumption~\ref{asm:step-size}-\eqref{asm:4-a2/b} we can conclude
\begin{align}
    &\lim_{T\rightarrow \infty} \EE{\sum_{t=1}^T \alpha(t) \dl(t)} \nonumber\\
    & = \lim_{T\rightarrow \infty} \sum_{t=1}^T \alpha(t) \EE{ \dl(t)} \nonumber\\
    & \leq  \sqrt{\frac{\cThree \cfone}{ \lambda}} \sum_{t=1}^\infty \alpha(t) \beta^{\frac{1}{2}}(t)
    + \sqrt{ \frac{2 \cFour \cftwo}{\lambda} } \sum_{t=1}^\infty  \frac{\alpha^2(t)}{\beta(t)} <\infty.
\end{align}
Using Monotone Convergence Theorem,
we have
\begin{align*}
    &\EE{\sum_{t=1}^\infty \alpha(t) \dl(t)}  = \lim_{T\rightarrow \infty} \EE{\sum_{t=1}^T \alpha(t) \dl(t)} < \infty,
\end{align*}
which implies 
\begin{align}\label{eq:sum-alpha-delta}
 \sum_{t=1}^\infty \alpha(t) \dl(t) < \infty,   
\end{align} 
almost surely. 

Now, recall random variables $v_1,v_2,\dots, v_B$ defined in  Section~\ref{sec:proof:variance}. We aim to show that these random variables are all zero, almost surely. We prove this claim by contradiction. Assume there exists some $\tau\in\{1,2,\dots, B\}$ such that  $p:=\Pr[v_\tau >0]>0$. Hence, by the continuity of measure, there  exists some (deterministic) $\epsilon>0$ such that $\Pr[v_\tau >\epsilon]>p/2>0$. Consider the event $\mathcal{A}=\{v_\tau >\epsilon\}$. Then $ \lim_{k\rightarrow\infty} \dl^2(\tau+kB) = v_\tau $ and $\dl(t) \geq 0$ imply that for all $\omega\in \mathcal{A}$, there exists some $\kz$ (possibly depending on~$\omega$) such that ${\dl(\tau+ k B ) \geq \sqrt{ \epsilon/2}}$ for $k\geq \kz$. 
Then, since $\{\alpha(t)\}$ is a non-increasing sequence, we have 
\begin{align}\label{eq:sum-al}
    \sum_{t=1}^\infty \alpha(t) \dl(t)  &\geq \sum_{k=0}^{\infty} \alpha(\tau+ k B) \dl(\tau+ k B)\nonumber\\
    & \geq \sum_{k=\kz}^{\infty} \alpha(\tau+ k B) \dl(\tau+ k B)\nonumber\\ 
    & \geq \sqrt{ \frac{\epsilon}{2}}\sum_{k=\kz}^{\infty} \alpha(\tau+ k B) \nonumber\\
    & \geq  \sqrt{ \frac{\epsilon}{2}} \sum_{k=\kz}^\infty 
     \frac{1}{B} \sum_{j= 0}^{B-1} \alpha(\tau+ kB+j)\nonumber\\
    & = \frac{1}{B}\sqrt{ \frac{\epsilon}{2}} \sum_{\ell= \tau+\kz B }^{\infty} \alpha(\ell)  = \infty,
\end{align}
where the last equality follows from Assumption~\ref{asm:step-size}-\eqref{asm:4-a}.
This implies that 
\[\Pr\left[\sum_{t=1}^\infty \alpha(t) \dl(t)=\infty\right]\geq \Pr(\mathcal{A})>\frac{p}{2}>0,\] which is in  contradiction with~\eqref{eq:sum-alpha-delta}. Therefore, we have
${v_1\!=\!v_2\!=\cdots=\!v_B\!=\!0}$, with probability $1$. 
\subsection{Average State Distance to an Optimal Point}
\label{sec:ave-to-opt}
Now, we derive an upper bound for the expected distance between the (weighted) average of the agents' states, i.e., ${\mx(t) = \bfr^T X(t) }$ and an arbitrary minimizer of $\bfx^\star\in\mathcal{X}^\star$ of the function $f(\bfx)$. Recall that  ${\bfr^T A(t) = \bfr^T}$. Hence, multiplying both sides of~\eqref{eqn:UpdateLawMatrixForm} by $\bfr^T$, subtracting $\bfx^\star$, and taking expectation, we arrive at
\begin{align}
\!\!\EE{\!\norm{\mx(t\!+\!1) \!-\! \bfx^\star}^2\! \md\F_t\!}& \!\!= \!\EE{\norm{\mx(t) \!+\! \bfr^TU(t) - \bfx^\star}^2\md\F_t} \nonumber \\
    & = \!\norm{\mx(t) \!-\! \bfx^\star}^2
\!+\!\EE{\norm{\bfr^TU(t)}^2\md\F_t}\!+\!2\!\left\langle\EE{\bfr^T U(t)\md\F_t} \!, \mx(t) \!-\! \bfx^\star\! \right\rangle\!.\label{eqn:ConsensusInit}
\end{align}
Using Lemma~\ref{lem:CauchyExtn} with $\theta=1$, we can bound the second term in~\eqref{eqn:ConsensusInit} as
\begin{align}\label{eq:b_a_Ft}
    \EE{\norm{\bfr^TU(t)}^2\md\F_t} &= \EE{\norm{\beta(t)\bfr^T E(t) - \alpha(t)\bfr^T \G(t)}^2\md\F_t}\nonumber\\
    & \leq\! 2\beta^2(t)\EE{\norm{\bfr^TE(t)}^2\md\F_t} \!+ \! 2\alpha^2(t)\! \norm{\bfr^T \G(t)}^2\!\!.
\end{align}
Note that $\ER{t}E^T(t)$ is an $n\times n$ matrix, and its $(i,j)$th entry is $\nt_i(t)\nt_j^T(t)$, which can be bounded using  Assumption~\ref{Assumption:Noise} as
\begin{align*}
       \lb\EE{|\ER{t}E^T(t)| \md\F_t}\rb_{ij}&=\EE{|\nt_i(t) \nt_j^T(t) | \md \F_t} \cr 
        &\leq \sqrt{\EE{\|\nt_i(t) \|^2 \md \F_t} \EE{\|\nt_j(t)\|^2 \md \F_t} }\cr
        &\leq \sm,
   \end{align*}
   for all $1\leq i,j\leq n$. Since $\bfr$ is a non-negative vector and $\bfr^T\ones =1$, for the first term in~\eqref{eq:b_a_Ft}, we have
  \begin{align}\label{eq:quant-noise}
        \EE{\left\|\bfr^T \ER{t}
        \right\|^2 \md\F_t}
         &=\bfr^T\EE{\ER{t}E^T(t)
         \md\F_t}  \bfr
        \nonumber\\
        &\leq\bfr^T\EE{ |\ER{t}E^T(t)|
         \md\F_t}  \bfr \nonumber\\
         &\leq \bfr^T (\sm \ones \ones^T) \bfr = \sm.
   \end{align}
   Similarly, from Assumption~\ref{Assumption:Function}-(c), we arrive at
   \begin{align}
        \left\|\bfr^T \G(t)
        \right\|^2 & = \bfr^T \G(t) \lb \G(t) \rb^T \bfr
        \nonumber\\
        &\leq \bfr^T (L^2 \ones \ones^T) \bfr = L^2.
        \label{eq:bnd-grad}
   \end{align}
Plugging~\eqref{eq:quant-noise} and~\eqref{eq:bnd-grad} into~\eqref{eq:b_a_Ft}, we can write
\begin{align}
    \EE{\norm{\bfr^T U(t)}_2^2 \md \F_t} &\leq 
     2\beta^2(t)\gamma + 2\alpha^2(t)L^2.\label{eqn:SecondTermBounded}
\end{align}
Recall that Assumption~\ref{Assumption:Noise} implies $ \EE{ \beta(t)\bfr^T\ER{t}\md\F_t} = 0 $. Using this fact and linearity of inner product, we can bound  the last term in~\eqref{eqn:ConsensusInit} as
\begin{align}
    \langle &\EE{\bfr^TU(t)\md\F_t},\mx(t) - \bfx^\star \rangle \nonumber\\
    & = \!\langle \EE{\beta(t)\bfr^T\! \ER{t}\md\F_t} \!-\!  \EE{\alpha(t)\bfr^T \G(t) \md\F_t}\!,\mx(t) \!-\! \bfx^\star \rangle\nonumber\\
    &= -\alpha(t) \langle\bfr^T \G(t) ,\mx(t) - \bfx^\star \rangle\nonumber\\
    & = - \alpha(t)\left\langle
    \sum\nolimits_{i=1}^n r_i \bfg_i(\bfx_i(t)) ,\mx(t) - \bfx^\star\right\rangle\nonumber\\
   & = -\alpha(t) \sum\nolimits_{i=1}^n\! r_i\langle \bfg_i(\bfx_i(t)) ,\mx(t) -\bfx^\star\rangle. \label{eq:grad_f_i}
\end{align} 
Let us consider each summand in~\eqref{eq:grad_f_i}, where we can write
\begin{align}\label{eq:grad_f_i_2}
    \langle \bfg_i(\bfx_i(t)) ,\mx(t) \!-\! \bfx^\star\rangle & = \langle \bfg_i(\bfx_i(t)),\mx(t) - \bfx_i(t)\rangle + \!\langle \bfg_i(\bfx_i(t)) ,\bfx_i(t) - \bfx^\star\rangle.
\end{align}
Using the Cauchy-Schwarz inequality and {Assumption~\ref{Assumption:Function}-(c)}, the first term in~\eqref{eq:grad_f_i_2} can be lower bounded as
\begin{align}\label{eq:grad_f_i_3}
     \langle \bfg_i(\bfx_i(t)),\mx(t) - \bfx_i(t)\rangle & \geq 
    -\| \bfg_i(\bfx_i(t)) \|\norm{\mx(t) - \bfx_i(t)}\nonumber\\
    & \geq 
    -L\norm{\mx(t) - \bfx_i(t)}.
\end{align}
From the convexity of $f_i(\cdot)$ in Assumption~\ref{Assumption:Function}-(a), for the second term in~\eqref{eq:grad_f_i_2} we have
\begin{align}\label{eq:grad_f_i_4}
    \langle \bfg_i(\bfx_i(t)),\bfx_i(t) - \bfx^\star\rangle & \geq f_i(\bfx_i(t)) - f_i(\bfx^\star)\nonumber\\
    &  = f_i(\mx(t))- f_i(\bfx^\star) + f_i(\bfx_i(t)) - f_i(\mx(t))  \nonumber\\
    & \geq f_i(\mx(t)) - f_i(\bfx^\star)\ - L\norm{\bfx_i(t) - \mx(t)},
\end{align}
where the last inequality follows from Assumption~\ref{Assumption:Function}-(c).
Therefore, substituting~\eqref{eq:grad_f_i_3} and~\eqref{eq:grad_f_i_4} into~\eqref{eq:grad_f_i_2}, we get
 \begin{align}\label{eq:inner-gr-i} 
    \langle \bfg_i(\bfx_i(t)) ,\mx(t) - \bfx^\star\rangle \geq &
    -2L\norm{\bfx_i(t) -\mx(t) } \nonumber\\
    & + f_i(\mx(t)) - f_i(\bfx^\star).
\end{align}
Replacing~\eqref{eq:inner-gr-i} in~\eqref{eq:grad_f_i}, and using the Cauchy-Schwarz inequality and the fact that $\sum_{i=1}^n r_i=1$, we have 
\begin{align}\label{eqn:ThirdTermFinal}
&\langle\EE{\bfr^TU(t)\md\F_t}\!,\mx(t)\! -\! \bfx^\star \rangle \nonumber\\
& \leq  2\alpha(t) L \! \sum_{i=1}^{n}r_i\norm{ \bfx_i(t) \!-\!\mx(t)} \!-\! \alpha(t) \! \sum_{i=1}^{n}r_i (f_i(\mx(t))\! -\!  f_i(\bfx^\star))\nonumber\\
&\leq  2\alpha(t) L  \sqrt{\sum_{i=1}^{n}r_i \cdot \sum_{i=1}^{n}r_i\norm{ \bfx_i(t) \!-\!\mx(t)}^2 }\nonumber\\
& \phantom{\leq }- \alpha(t)  \sum_{i=1}^{n}r_i (f_i(\mx(t))\! -\!  f_i(\bfx^\star))\nonumber\\
&=  2\alpha(t) L  \dl(t) - \alpha(t) \left( f(\mx(t))\! -\! f(\bfx^\star)\right).
\end{align}
Plugging~\eqref{eqn:SecondTermBounded} and~\eqref{eqn:ThirdTermFinal} into~\eqref{eqn:ConsensusInit}, we get
\begin{align*}
    \mathbb{E}\left[\hspace{-1pt}\norm{\mx(t\!+\!1) \!-\! \bfx^\star}^2\md\F_t\right]
    &\hspace{-2pt} \leq\! \norm{\mx(t)\! - \!\bfx^\star}^2 \!+\! 2(\beta^2(t)\gamma \!+\! \alpha^2(t)L^2)\\
       &\phantom{\leq} + 4\alpha(t)L\dl(t) 
        \hspace{-2pt}-\!2\alpha(t)(f(\mx(t)) \!-\! f(\bfx^\star)),
\end{align*}
which is identical to the inequality in  Theorem~\ref{thm:a_s_conv}, with $\bfy(t) = \mx(t)$, $\zeta(t) = 0$, $\xi(t) = 2\alpha(t)$, and
\begin{align*}
     z(t) = 4L\alpha(t)\dl(t) + 2\gamma  \beta^2(t)\!+\! 2L^2\alpha^2(t).
\end{align*}
Note that $\zeta(t)$,  $\xi(t)$, and and $z(t)$ are all none-negative for $t\geq 1$, and $\sum_{t=1}^{\infty}\zeta(t)=0< \infty$. Moreover,  Assumption~\ref{asm:step-size}-\eqref{asm:4-a}
 implies $\sum_{t=1}^{\infty}\xi(t) = 2 \sum_{t=1}^{\infty} \alpha(t) =\infty$. Finally, from~\eqref{eq:sum-alpha-delta} and Assumption~\ref{asm:step-size}-\ref{asm:4-a2-b2} we have $\sum_{t=1}^{\infty} z(t)<\infty$. Therefore, we can apply  Theorem~\ref{thm:a_s_conv}, and conclude that $\{\mx(t)\}$ converges to some $\tilde{\bfx} \in \mathcal{X}^\star$, almost surely.  
\subsection{Almost Sure Convergence of State to an Optimal Point}
\label{sec:state-to-average}
In Section~\ref{sec:var} we showed that under the Assumptions~\ref{Assumption:Noise}-\ref{asm:step-size}, we have $\lim_{t\rightarrow\infty} \dl(t) = \lim_{t\rightarrow\infty} \Nr{X(t) - \ones \mx(t)} =0$. This implies $ \lim_{t\rightarrow \infty} \Nr{\bfx_i(t) -\mx(t)} = 0$ for every $i\in [n]$. Moreover, we proved that $\mx(t)$ converges to some $\tilde{\bfx} \in \mathcal{X}^\star$, almost surely, in Section~\ref{sec:ave-to-opt}. Combining these two results, immediately conclude the claim of Theorem~\ref{thm:almost_sure}. 
\section{Proof of Proposition~\ref{prop:mu_nu}}\label{sec:proof-prop}
In this section, we prove Proposition~\ref{prop:mu_nu}. We only need to show that step-sizes ${\beta(t) = \frac{\beta_0}{t^\mu}}$ and ${\alpha(t) = \frac{\alpha_0}{t^\nu}}$ with ${\frac{1}{2}<\mu \leq 1}$ and ${\frac{1}{2}(1\!+\!\mu)< \nu  \leq 1}$ satisfy all conditions in Assumption~\ref{asm:step-size}.
First, from $\bzr\leq 1$, we get $\beta(t)= \frac{\bzr}{t^\mu}\leq 1$ for all $t\geq 1$. Using Lemma~\ref{lm:sum-exp} with $\nu \leq 1$, we have $\sum_{t=1}^{\infty} \alpha(t)=\sum_{t=1}^{\infty}\frac{\azr}{t^\nu} =\infty$.
Moreover, from Lemma~\ref{lm:sum-exp} with $\mu>\frac{1}{2}$, $\nu>\frac{1}{2}$, and $2\nu-\mu>1$ we arrive at
\begin{align*}
   & \sum_{t=1}^{\infty} \beta^2(t) = \sum_{t=1}^{\infty} \frac{\bzr^2}{t^{2\mu}}< \infty,\\
   & \sum_{t=1}^{\infty}\alpha^2(t) = \sum_{t=1}^{\infty} \frac{\azr^2}{t^{2\nu}} < \infty,\\
   & \sum_{t=1}^{\infty}\frac{\alpha^2(t)}{\beta(t)} = \sum_{t=1}^{\infty} \frac{\azr^2}{\bzr}\frac{1}{t^{2\nu-\mu}} < \infty.
\end{align*}
Now, we need to show that for some positive constants $\cOne\!<\!\frac{\lambda}{2}$, $\cTwo\!<\!\frac{\lambda}{4}$, and $\tz \geq 1$ we have $-\Delta\beta(t)\leq \cOne \beta^2(t)$ and $-\Delta\alpha(t)\leq \cTwo \alpha(t)\beta(t)$ for all $t\geq \tz$. 

Using the mean value theorem for the function $\beta(t) = \frac{\bzr}{t^\mu}$ we have ${\beta(t+1) - \beta(t) = \beta'(\zeta_1)}$ for some $\zeta_1 \in [t, t+1]$. Therefore, we arrive at
\begin{align*}
   - \Delta \beta(t) & = \beta(t) - \beta(t+1)\\
   & = -\beta'(\zeta_1) = \mu\frac{\bzr}{\zeta_1^{\mu+1}} \leq \mu \frac{\bzr}{t^{\mu+1}} \leq \cOne \frac{\bzr^2}{ t^{2\mu}} = \cOne \beta^2(t),
\end{align*}
where the latter holds for $t \geq \tone:=\lp \frac{\mu}{ \bzr \cOne}\rp^{\frac{1}{1-\mu}}$ provided that $\mu< 1$. Similarly, for the the function $\alpha(t)=\frac{\azr}{t^{\nu}}$ we get $\alpha(t+1)-\alpha(t) = \alpha'(\zeta_2)$ for some $\zeta_2\in [t,t+1]$. Hence, we can write
\begin{align*}
   - \Delta \alpha(t) & = \alpha(t) - \alpha(t+1)\\
   & =\! -\alpha'(\zeta_2) = \nu\frac{\azr}{\zeta_2^{\nu+1}} \!\leq\! \nu \frac{\azr}{t^{\nu+1}} \!\leq\! \cTwo \frac{\azr\bzr}{t^{\nu+\mu}} \!=\! \cTwo \alpha(t)\beta(t),
\end{align*}
for every $t\geq \ttwo$ where $\ttwo :=  \lp \frac{\nu}{ \bzr \cTwo}\rp^{\frac{1}{1-\mu}}$. Therefore, for any pair of (fixed) positive constants $\cOne\!<\!\frac{\lambda}{2}$ and $\cTwo\!<\!\frac{\lambda}{4}$ we have $-\Delta\beta(t)\leq \cOne \beta^2(t)$ and $-\Delta\alpha(t)\leq \cTwo \alpha(t)\beta(t)$ for all ${t\geq \tz=\max \lp \tone,\ttwo\rp}$. This shows that Assumption~\ref{asm:step-size}-\eqref{asm:4-Db}-\eqref{asm:4-Da} are satisfied for sufficiently large $t$, regardless of the dynamic parameters. This completes the proof of Proposition~\ref{prop:mu_nu}.

\section{Conclusion} \label{sec:conclusion}
In this work, we have studied distributed optimization over time-varying networks suffering from noisy/lossy communication between the agents using a two-time-scale consensus-based algorithm. We have identified sufficient conditions for general step-sizes sequences for the two-time-scales, including damping mechanisms for the imperfect received information from neighboring agents as well as the local loss functions' gradients, to guarantee the algorithm's \textit{almost sure convergence} for \textit{convex} cost functions. Furthermore, we used this result to characterize conditions on practical step-size sequences that enables almost sure convergence in this setting. 

Future efforts in this area may include identifying \textit{necessary} conditions on the step-size sequences for the almost sure convergence of the algorithm. In particular, it is interesting to study the discrepancy between the two converging regions for convex and strongly convex settings in Fig.~\ref{fig:Lu}. In addition, the extension of this work to distributed online learning algorithms is of future interest. 



\bibliographystyle{IEEEtran}
\bibliography{ref}

\appendix
\section*{Appendix: Proof of Preliminaries}
{\it Proof of Lemma~\ref{lemma:transition2}}:
	Due to the separable nature of $\Nr{\cdot}$, i.e.,  $\Nr{U}^2=\sum_{j=1}^d\Nr{U^{j}}^2$, without loss of generality, we may assume that $d=1$. Thus, let $U=\bfu\in \R^n$. Define $\V:\R^n\to \R^+$ by 
	    \begin{align}
	        \V(\bfu) & := \Nr{(I-\ones \bfr^T)\bfu}^2\nonumber\\
	        & = \Nr{\bfu-\ones \bfr^T\bfu}^2=\sum_{i=1}^n r_i(u_i-\bfr^T\bfu)^2.
	    \end{align}
	    Let us denote ${\bfu(s)\!=\!\bfu=\begin{bmatrix} u_1 & u_2 & \cdots & u_n \end{bmatrix}}$ and  $\bfu(k+1)= A(k+1)\bfu(k)$. In addition with a slight abuse of notation, we  denote $\V(\bfu(k))$ by $\V(k)$ for $k=s,\ldots,t$. 
	    
	    Using Theorem 1 in  \cite{touri2011existence}, we have 
	    \begin{align}\label{eqn:decreaseV}
	       \V(t)=\V(s)-\sum_{k=s+1}^{t}\sum_{i<j}H_{ij}(k)(u_i(k)-u_j(k))^2,
	    \end{align}
	    for $t> s$, where $H(k)=A^T(k)\diag(\bfr)A(k)$, is a  non-negative matrix.  Then, setting $\bfu(s) = \bfu $ and ${\bfu(t-1) = \Phi(t:s) \bfu(s)}$ for $t>s$, we have
	    \begin{align}\label{eq:prf:lm-tran-1}
	        \Nr{(\Phi(t\hspace{-1pt}:\hspace{-1pt}s)-\ones \bfr^T) \bfu)}^2 &\stackrel{\rm{(a)}}{=} 
	        \Nr{(I-\ones \bfr^T) \Phi(t\hspace{-1pt}:\hspace{-1pt}s)\bfu)}^2\nonumber\\
	        & = \Nr{(I-\ones \bfr^T) \bfu(t-\hspace{-1pt}1)}^2\nonumber\\
	        &= \V(t-1) \leq \V(s) \nonumber\\
	        & = \Nr{\bfu-\ones \bfr^T \bfu)}^2 \leq \Nr{\bfu}^2,
	    \end{align}
	    where {\rm (a)} follows from  Assumption~\ref{Assumption:Graph}-(a) and the fact that
     \begin{align}\label{eq:Ak}
         {A(k) = (1-\beta(k))I + \beta(k) W(k)}
     \end{align} which imply ${\bfr^T \Phi(t:s) =\bfr^T}$, and the inequality in {\rm (b)} is due to the fact that  \[\Nr{\bfu - \ones \bfr^T \bfu}^2 + \Nr{\ones \bfr^T \bfu}^2 = \Nr{\bfu}^2.\] This implies the first claim of the lemma in~\ref{eq:lm-tran-1}.
	    
	    Furthermore, since $A(k)$ is a non-negative matrix, we have ${H(k)\geq \rmin A^T(k)A(k)}$, for $k=s+1,\ldots,t$. Also, since $A(k)$ satisfies \eqref{eq:Ak}, then  Assumption~\ref{Assumption:Graph}-(b) implies that the minimum non-zero elements of $A(k)$ are bounded bellow by $\minW\beta(k)$. Therefore, since $\beta(k)$ is non-increasing, on the window $k=s+1,\ldots,s+B$, the minimum non-zero elements of $A(k)$  for $k$ in this window are lower bounded by $\minW\beta(s+B)$. Without loss of generality, assume that the entries of $\bfu$ are sorted, i.e., $u_1\leq \ldots\leq u_n$, otherwise, we can relabel the agents (rows and columns of $A(k)$s and $\bfu$ to achieve this).  Therefore, by Lemma~8 in \cite{nedic2008distributed}, for \eqref{eqn:decreaseV}, we have 
	    \begin{align}\label{eqn:decreaseV2}
	       \V(s+B) &\leq\V(s)\!-\! \rmin\sum_{k=s+1}^{s+B}\sum_{i<j}[A^T(k)A(k)]_{ij}(u_i(k)-u_j(k))^2\nonumber\\
	       &\leq \V(s)-\frac{\minW\rmin}{2}\beta(s+B)\sum_{\ell=1}^{n-1}(u_{\ell+1}-u_{\ell})^2.
	    \end{align}
	    We may comment here that although Lemma~8 in \cite{nedic2008distributed} is written for doubly stochastic matrices, and its statement is about the decrease of $\V(\bfx)$ for the special case of  $\bfr=\frac{1}{n}\ones$, but in fact, at the core of its proof, it is a result on bounding 
     \[\sum_{k=s+1}^{s+B}\sum_{i<j}[A^T(k)A(k)]_{ij}(u_i(k)-u_j(k))^2\] for a sequence of {$B$-connected} stochastic matrices $A(k)$ in terms of the minimum non-zero entries of stochastic matrices ${A(s+1),\ldots,A(s+B)}$. 
	    
	    Next, we will show that $\sum_{\ell=1}^{n-1}(u_{\ell+1}-u_{\ell})^2\geq n^{-2}\V(\bfu)$. This argument adapts a similar argument used in the proof of Theorem~18 in \cite{nedic2008distributed} to the general $\V(\cdot)$. 
	    
	    For a $\bfv\in \R^n$ with $\V(\bfv)>0$, define the quotient 
	    \begin{align}
	        h(\bfv)=\frac{\sum_{\ell=1}^{n-1}(v_{\ell+1}-v_{\ell})^2}{\sum_{i=1}^nr_i(v_i-\bfr^T\bfv)^2}=\frac{\sum_{\ell=1}^{n-1}(v_{\ell+1}-v_{\ell})^2}{\V(\bfv)}. 
	    \end{align}
	    Note that $h(\bfv)$ is invariant under scaling and translations by all-one vector, i.e., $h(\omega \bfv )=h(\bfv)$ for all non-zero $\omega\in \R$ and $h(\bfv+\omega \ones)=h(\bfv)$ for all $\omega\in \R$. Therefore, 
	    \begin{align}\label{eqn:hv}
	        \min_{\substack{v_1\leq v_2\leq \cdots\leq v_n\\\V(\bfv)\not=0}}h(\bfv)&=\min_{\substack{v_1\leq v_2\leq \cdots\leq v_n\\ \bfr^T\bfv=0,\V(\bfv)=1}}h(\bfv)\cr 
	        &=\min_{\substack{v_1\leq  v_2\leq \cdots\leq v_n\\\bfr^T\bfv=0,\V(\bfv)=1}}\sum_{\ell=1}^{n-1}(v_{\ell+1}-v_{\ell})^2. 
	    \end{align}
	    
	    Since $\bfr$ is a stochastic vector, then for a vector $\bfv$ with  ${v_1\leq \ldots\leq v_n}$ and $\bfr^T\bfv=0$, we would have ${v_1\leq \bfr^T\bfv=0\leq v_n}$. On the other hand, the fact that ${\V(\bfv)=\sum_{i=1}^nr_iv^2_i=1}$ would imply ${\max(|v_1|,|v_n|)\geq \frac{1}{\sqrt{n}}}$. Let us consider the difference sequence  $\hat{v}_\ell=v_{\ell+1}-v_{\ell}$  for $\ell=1,\ldots, n-1$, for which we have $\sum_{i=1}^{n-1} \hat{v}_\ell = v_n - v_1 \geq v_n \geq \frac{1}{\sqrt{n}}$. Therefore, the optimization problem \eqref{eqn:hv} can be  rewritten as
        \begin{align}\label{eqn:hv2}
	        \min_{\substack{v_1\leq v_2\leq \cdots\leq v_n\\\V(\bfv)\not=0}}h(\bfv)&=\min_{\substack{v_1\leq v_2\leq \cdots\leq v_n\\\bfr^T\bfv=0,\V(\bfv)=1}}\sum_{\ell=1}^{n-1}(v_{\ell+1}-v_{\ell})^2\cr 
	        &\geq 
	        \min_{\substack{\hat{v}_1,\ldots,\hat{v}_{n-1}\geq 0\\\sum_{i=1}^{n-1}\hat{v}_i\geq \frac{1}{\sqrt{n}}}}\sum_{\ell=1}^{n-1}\hat{v}_\ell^2. 
	    \end{align}
	    Using the Cauchy-Schwarz inequality, we get 
	    ${\lp \sum_{\ell=1}^{n-1} \hat{v}_\ell^2 \rp \cdot \lp \sum_{\ell=1}^{n-1} 1^2\rp \geq \big( \sum_{\ell=1}^{n-1} \hat{v}_\ell \big)^2 \geq \big( \frac{1}{\sqrt{n}} \big)^2 = \frac{1}{n}}$. Hence, 
	    
	    \begin{align}
	         \min_{\substack{v_1\leq v_2\leq \ldots\leq v_n\\
	         \V(\bfv)\not=0}}h(\bfv)\geq \frac{1}{n(n-1)}\geq \frac{1}{n^2}.
	    \end{align}
        Thus, for ${v_1 \!\leq\! \ldots \!\leq\! v_n}$, we have $\sum_{\ell=1}^{n-1}(v_{\ell+1}\!-\!v_{\ell})^2\!\geq \! n^{-2}\V(\bfv)$ (note that this inequality also holds for $\bfv\in\R^n$ with $\V(\bfv)=0$). 
	    Using this fact in \eqref{eqn:decreaseV2} implies  
	    \begin{align}\label{eqn:Vdot}
	       \V(s+B) \leq \lp 1-\frac{\minW\rmin}{2n^2}\beta(s+B)\rp\V(s).
	    \end{align}
       Therefore, similar to~\eqref{eq:prf:lm-tran-1}, we can continue  from~\eqref{eqn:Vdot} and write
	   \begin{align*}
	       \big\|\!\big( \Phi(s\!+\!B\!+\!1\!:\!s) \!-\!\ones \bfr^T\big) \bfu \big\|_\bfr^2 & \!=\! 
	       \big\| \big( I\!-\! \ones \bfr^T\big) \Phi(s\!+\!B\!+\!1\!:\!s) \bfu \big\|_\bfr^2\\
	       &=\Nr{(I - \ones \bfr^T) \bfu(s+B)}^2\\
	       & = \V(s+B) \\
	       & \leq  \lp 1-\lambda B\beta(s\!+\!B)\rp\V(s)\\
	       & = \!\lp 1\!-\!\lambda B\beta(s\!+\!B)\rp\!\Nr{\bfu \!-\! \ones \bfr^T\! \bfu}^2 \\
	       & \leq  \lp 1-\lambda B\beta(s+B)\rp\Nr{\bfu}^2. 
	   \end{align*}
	   Applying this inequality on each column of a matrix $U$, we can conclude the same result for matrices. This completes the proof of the lemma. \hfill $\blacksquare$

    {\it Proof of Lemma~\ref{lm:sum_pro_g}}: 
    We first define a sequence $g(t)$ via
    \begin{align}\label{eq:lm:sum-prod:rec}
        g(t+1)= (1-q(t)) g(t) + p(t),
    \end{align}
    for $t\geq 1$ and $g(1)=0$. Then, we can verify that 
    \begin{align*}
        g(t) = \sum_{s=1}^{t-1} \lb p(s) \prod_{k=s+1}^{t-1} (1-q(k))\rb.
    \end{align*}
    for all values of $t\geq 1$. We set $S=\max\lc\frac{g(t_0) q(t_0)}{p(t_0)}, \frac{1}{1-A}\rc$, and we aim to show that $g(t) \leq S p(t) / q(t)$  for every $t\geq t_0$. 
    
    We use induction to prove the claim. First note that for $t=t_0$, we have $g(t_0) \leq S p(t_0)/q(t_0)$. Assume the claim holds for $t$. Then, for $t+1$ we have 
    \begin{align*}
        \sum_{s=1}^{t} \lb p(s) \prod_{k=s+1}^{t} (1-q(k))\rb & = g(t+1)\\
        & = (1-q(t)) g(t) + p(t) \\
        & \leq (1-q(t)) S\frac{p(t)}{q(t)} + p(t)\\
        & = S\frac{p(t)}{q(t)} - (S-1) p(t). 
    \end{align*}
    Thus, in order to show that $g(t+1) \leq S \frac{p(t+1)}{q(t+1)}$, it suffice to show that 
    \begin{align*}
         \frac{p(t)}{q(t)} -\frac{p(t+1)}{q(t+1)} \leq \frac{S-1}{S} p(t).
    \end{align*}
    To this end, we can write 
    \begin{align*}
         \frac{p(t)}{q(t)} -\frac{p(t+1)}{q(t+1)}
         &= \frac{p(t)}{q(t)} - \frac{p(t+1)}{q(t)} +\frac{p(t+1)}{q(t)} -\frac{p(t+1)}{q(t+1)} \\
         &=\frac{-\Delta p(t)}{q(t)} + p(t+1) \frac{q(t+1)- q(t)}{q(t) q(t+1)} \\
         &\leq \frac{-\Delta p(t)}{q(t)}\\
         &\leq \frac{A p(t) q(t)}{q(t)} \cr
         &\leq \frac{S-1}{S} p(t),
    \end{align*}
    where the first inequality holds since $\{q(t)\}$ is a non-increasing sequence, the second inequality follows from~\eqref{eq:sum_prod_con}, and the last inequality holds since $S\geq \frac{1}{1-A}$. This completes the proof of the lemma. \hfill $\blacksquare$ 

\end{document}